\documentclass[11pt]{article}
\usepackage[matrix,arrow,curve,cmtip]{xy}
\usepackage{amssymb}
\usepackage{latexsym}
\usepackage{theorem}
\usepackage[a4paper,width=160mm,height=230mm]{geometry}
\usepackage{titlesec}
\usepackage{graphicx}
\usepackage[applemac]{inputenc}
\usepackage{eucal}

\titleformat{\section}[hang]%
{\bf\Large}{\thesection.}{1ex}{}%

\titleformat{\subsection}[hang]%
{\bfseries\normalsize}{\thesubsection.}{1ex}{}

\def\to{\mbox{$\xymatrix@1@C=5mm{\ar@{->}[r]&}$}}
\def\halfcirc{\begin{picture}(0,0)\put(0,3){\oval(4,4)[l]}\end{picture}}
\def\incl{\mbox{$\xymatrix@1@C=5mm{\ar@{->}[r]|<{\halfcirc}&}$}}
\def\tto{\mbox{$\xymatrix@1@C=5mm{\ar@{=>}[r]&}$}}
\def\distsign{\begin{picture}(0,0)\put(0,0){\circle{4}}\end{picture}}
\def\dist{\mbox{$\xymatrix@1@C=5mm{\ar@{->}[r]|{\distsign}&}$}}
\def\spansign{\begin{picture}(0,0)\put(0,-3){\line(0,1){6}}\end{picture}}
\def\span{\mbox{$\xymatrix@1@C=5mm{\ar@{->}[r]|{\spansign}&}$}}
\def\criblesign{\begin{picture}(0,0)\put(-1,-3){\line(0,1){6}}\put(1,-3){\line(0,1){6}}\end{picture}}
\def\crible{\mbox{$\xymatrix@1@C=5mm{\ar@{->}[r]|{\criblesign}&}$}}

\def\inlineadj#1#2{\mbox{$\xymatrix@C=15mm{\ar@{}[r]|{\bot}\ar@<1mm>@/^2mm/[r]^{{#1}} &\ar@<1mm>@/^2mm/[l]^{{#2}}}$}}

\newtheorem{theorem}{Theorem}[section]
\newtheorem{lemma}[theorem]{Lemma}
\newtheorem{definition}[theorem]{Definition} 
\newtheorem{proposition}[theorem]{Proposition}
\newtheorem{corollary}[theorem]{Corollary}
{\theorembodyfont{\upshape}\newtheorem{example}[theorem]{Example}}
\newcommand{\proof}{\noindent {\em Proof\ }: }
\def\endofproof{$\mbox{ }\hfill\Box$\par\vspace{1.8mm}\noindent}

\def\Rel{{\sf Rel}}
\def\o{^{\sf o}}
\def\bigmid{~\Big|~}
\def\Sh{{\sf Sh}}
\def\si{_{\sf si}}

\def\ssi{_{\sf ssi}}

\def\Sym{{\sf Sym}}
\def\Matr{{\sf Matr}}
\def\Cat{{\sf Cat}}
\def\Mod{{\sf Mod}}
\def\Ord{{\sf Ord}}
\def\Idl{{\sf Idl}}
\def\:{\colon}

\def\2{{\bf 2}}
\def\Set{{\sf Sh}}
\def\op{^{\sf op}}
\def\dom{{\sf dom}}

\def\Sup{{\sf Sup}}
\def\Dist{{\sf Dist}}
\def\Map{{\sf Map}}

\def\cc{_{\sf cc}}
\def\sc{_{\sf sc}}

\def\s{_{\sf s}}
\def\P{{\cal P}}
\def\C{{\cal C}}
\def\E{{\cal E}}
\def\A{{\cal A}}

\def\K{{\cal K}}
\def\Q{{\cal Q}}
\def\R{{\cal R}}

\def\V{{\cal V}}
\def\O{{\cal O}}
\def\topE{{\cal T}}
\def\bbA{\mathbb{A}}
\def\bbB{\mathbb{B}}
\def\bbC{\mathbb{C}}

\def\bbR{\mathbb{R}}

\def\bbP{\mathbb{P}}
\def\bbQ{\mathbb{Q}}
\def\tensor{\otimes}
\def\<{\langle}
\def\>{\rangle}

\def\Hilb{{\sf Hilb}}
\def\Proj{{\sf ProjMatr}}
\def\eqref#1{(\ref{#1})}

\def\m{^{\sf m}}

\title{Grothendieck quantaloids for allegories of enriched categories}
\author{Hans Heymans\footnote{Department of Mathematics and Computer Science, University of Antwerp, Middelheimlaan 1, 2020 Antwerpen, Belgium, {\tt hans.heymans@ua.ac.be}}\ \ and Isar Stubbe\footnote{Laboratoire de Math\'ematiques Pures et Appliqu\'ees, Universit\'e du Littoral-C\^ote d'Opale, 50 rue F. Buisson, 62228 Calais, France, {\tt isar.stubbe@lmpa.univ-littoral.fr}}}
\date{Submitted December 30, 2011; revised June 11, 2012; finalised June 26, 2012.}

\begin{document}

\maketitle

\begin{quote}{\small 
{\bf Abstract.} For any small involutive quantaloid $\Q$ we define, in terms of symmetric quantaloid-enriched categories, an involutive quantaloid $\Rel(\Q)$ of $\Q$-sheaves and relations, and a category $\Set(\Q)$ of $\Q$-sheaves and functions; the latter is equivalent to the category of symmetric maps in the former. We prove that $\Rel(\Q)$ is the category of relations in a topos if and only if $\Q$ is a modular, locally localic and weakly semi-simple quantaloid; in this case we call $\Q$ a Grothendieck quantaloid. It follows that $\Set(\Q)$ is a Grothendieck topos whenever $\Q$ is a Grothendieck quantaloid. Any locale $L$ is a Grothendieck quantale, and $\Set(L)$ is the topos of sheaves on $L$. Any small quantaloid of closed cribles is a Grothendieck quantaloid, and if $\Q$ is the quantaloid of closed cribles in a Grothendieck site $(\C,J)$ then $\Set(\Q)$ is equivalent to the topos $\Sh(\C,J)$. Any inverse quantal frame is a Grothendieck quantale, and if $\O(G)$ is the inverse quantal frame naturally associated with an étale groupoid $G$ then $\Set(\O(G))$ is the classifying topos of $G$.
}\end{quote}

\section{Introduction}\label{A}

A topos arising as the category of left adjoints in a locally ordered category, is the subject of P. Freyd and A. Scedov's [1990] study of {\em allegories}. More precisely, an allegory $\A$ is a modular locally ordered 2-category whose hom-posets have binary intersections; taking left adjoints (also known as ``maps'') in an allegory $\A$ thus produces a category $\Map(\A)$; and the interesting case is where the latter category is in fact a topos. Thus, in Freyd and Scedrov's own words, allegories ``are to binary relations between sets as categories are to functions between sets''. In practice, those interesting allegories arise most often as universal constructions on much smaller sub-allegories which are easier to describe explicitly. Freyd and Scedrov [1990] (but see also [Johnstone, 2002, A3]) give several theorems to this effect.

In [1982], R. Walters proved that any small site $(\C,J)$ gives rise to a small quantaloid $\R(\C,J)$ in such a way that the topos $\Sh(\C,J)$ is equivalent to the category of Cauchy-complete symmetric $\R(\C,J)$-enriched categories and functors. But the latter category is further equivalent to the category of all symmetric $\R(\C,J)$-categories and left adjoint distributors, and the quantaloid $\Sym\Dist(\R(\C,J))$ of all symmetric $\R(\C,J)$-enriched categories and all distributors is modular. In other words, the topos $\Sh(\C,J)$ is the category of maps in the allegory of symmetric $\R(\C,J)$-enriched categories and distributors--- which thus qualifies as an ``interesting'' allegory. 

In this paper we shall explain more precisely how ``sheaves via quantaloid-enrichment'' fit with ``toposes via allegories''. To that end, we define in Section \ref{C}, for any involutive quantaloid $\Q$, a new involutive quantaloid $\Rel(\Q)$, to be thought of as the locally posetal 2-category of ``$\Q$-sheaves and relations'', and a new category $\Set(\Q)$, to be thought of as the category of ``$\Q$-sheaves and functions''. The objects of $\Rel(\Q)$, resp.\ $\Set(\Q)$, are particular symmetric quantaloid-enriched categories, and the morphisms are distributors, resp.\ functors; the relation between the two is that $\Set(\Q)$ is the category of {\em symmetric} left adjoints in $\Rel(\Q)$. For appropriate $\Q$, these $\Q$-sheaves are, among the $\Q$-orders of [Stubbe, 2005b], precisely the {\em symmetric} ones.

We show in Section \ref{D} that, if $\Q=\R(\C,J)$, then $\Set(\Q)$ is equivalent to $\Sh(\C,J)$ and $\Rel(\Q)$ is equivalent to $\Rel(\Sh(\C,J))$; thus we recover and refine Walters' [1982] insight. More generally, we prove in Section \ref{E} that $\Rel(\Q)$ is equivalent to $\Rel(\topE)$ for some Grothendieck topos $\topE$ (and thus $\Set(\Q)$ is equivalent to $\topE$) if and only if $\Q$ is a modular, locally localic and weakly semi-simple quantaloid; we call these {\em Grothendieck quantaloids}. In other words, these Grothendieck quantaloids are precisely those for which the $\Q$-sheaves and relations form an ``interesting allegory''. Locales and inverse quantal frames [Resende, 2007, 2012] are examples of Grothendieck quantales. If $L$ is a locale, then $\Set(L)$ is in fact the topos of sheaves on $L$. And if $\O(G)$ is the inverse quantal frame associated to an étale groupoid $G$ [Resende, 2007], then it follows from [Heymans and Stubbe, 2009b; Resende, 2012] that $\Set(\O(G))$ is the classifying topos of that groupoid.

\section{Sheaves on an involutive quantaloid}\label{C}

The new notions that we will present at the end of this section draw heavily on the theory of quantaloid-enriched categories. For self-containedness we present some preliminaries in the first couple of subsections. For more details and for the many appropriate historical references we refer to [Stubbe, 2005a; Heymans and Stubbe, 2011].

\subsection*{Enrichment, involution and symmetry}

A {\em quantaloid} $\Q$ is, by definition, a category enriched in the symmetric monoidal closed category $\Sup$ of complete lattices and supremum-preserving functions; and a {\em homomorphism} $F\:\Q\to\R$ of quantaloids is a $\Sup$-enriched functor. An {\em involution} on a quantaloid $\Q$ is a homomorphism $(-)\o\:\Q\op\to\Q$ which is the identity on objects and satisfies $f^{\sf oo}=f$ for any morphism $f$ in $\Q$. The pair $(\Q,(-)\o)$ is then said to form an {\em involutive quantaloid}; we shall often simply speak of ``an involutive quantaloid $\Q$'', leaving the notation for the involution understood. When both $\Q$ and $\R$ are involutive quantaloids, then we say that $F\:\Q\to\R$ is a homomorphism of involutive quantaloids when it is a homomorphism such that $F(f\o)=(Ff)\o$. 

Whenever a morphism $f\:A\to B$ in a quantaloid (or in a locally ordered category, for that matter) is supposed to be a left adjoint, we write $f^*$ for its right adjoint. A {\em symmetric left adjoint} in an involutive quantaloid $\Q$ is a left adjoint whose right adjoint is its involute: $f^*=f\o$. Precisely as we write $\Map(\Q)$ for the category of left adjoints in $\Q$, we write $\Sym\Map(\Q)$ for the category of symmetric left adjoints. 

A {\em category $\bbA$ enriched in a quantaloid $\Q$} consists of a set $\bbA_0$ of objects, each $x\in\bbA_0$ having a type $ta\in\Q_0$, and for any $x,y\in\bbA_0$ there is a hom-arrow $\bbA(y,x)\:tx\to ty$ in $\Q$, subject to associativity and unit requirements: $\bbA(z,y)\circ\bbA(y,x)\leq\bbA(z,x)$ and $1_{tx}\leq\bbA(x,x)$ for all $x,y,z\in\bbA_0$. A {\em functor} $F\:\bbA\to\bbB$ between such $\Q$-categories is an object-map $x\mapsto Fx$ such that $tx=t(Fx)$ and $\bbA(y,x)\leq\bbB(Fy,Fx)$ for all $x,y\in\bbA$. Such a functor is smaller than a functor $G\:\bbA\to\bbB$ if $1_{tx}\leq\bbB(Fx,Gx)$ for every $x\in\bbA$. With obvious composition one gets a locally ordered 2-category $\Cat(\Q)$ of $\Q$-categories and functors.

For two objects $x,y\in\bbA$, the hom-arrows $\bbA(y,x)$ and $\bbA(x,y)$ go in opposite directions. Hence, to formulate a notion of ``symmetry'' for $\Q$-categories, it is far too strong to require $\bbA(y,x)=\bbA(x,y)$. Instead, at least for involutive quantaloids, a $\Q$-category $\bbA$ is {\em symmetric} when $\bbA(x,y)=\bbA(y,x)\o$ for every two objects $x,y\in\bbA$ [Betti and Walters, 1982]. We shall write $\Sym\Cat(\Q)$ for the full sub-2-category of $\Cat(\Q)$ determined by the symmetric $\Q$-categories (in which the local order is in fact symmetric, but not anti-symmetric).

A {\em distributor} $\Phi\:\bbA\dist\bbB$ between $\Q$-categories consists of arrows $\Phi(y,x)\:tx\to ty$ in $\Q$, one for each $(x,y)\in\bbA_0\times\bbB_0$, subject to two action requirements: $\bbB(y',y)\circ\Phi(y,x)\leq\Phi(y',x)$ and $\Phi(y,x)\circ\bbA(x,x')\leq\Phi(y,x')$ for all $y,y'\in\bbB_0$ and $x,x'\in\bbA_0$. The composite of such a distributor with another $\Psi\:\bbB\dist\bbC$ is written as $\Psi\tensor\Phi\:\bbA\dist\bbC$, and its elements are
$$(\Psi\tensor\Phi)(z,x)=\bigvee_{y\in\bbB_0}\Psi(z,y)\circ\Phi(y,x)$$
for $x\in\bbA_0$ and $z\in\bbC_0$. Parallel distributors can be compared elementwise, and in fact one gets a (large) quantaloid $\Dist(\Q)$ of $\Q$-categories and distributors. Each functor $F\:\bbA\to\bbB$ determines an adjoint pair of distributors: $\bbB(-,F-)\:\bbA\dist\bbB$, with elements $\bbB(y,Fx)$ for $(x,y)\in\bbA_0\times\bbB_0$, is left adjoint to $\bbB(F-,-)\:\bbB\dist\bbA$ in the quantaloid $\Dist(\Q)$. These distributors are said to be {\em represented by $F$}. More generally, a (necessarily left adjoint) distributor $\Phi\:\bbA\dist\bbB$ is {\em representable} if there exists a (necessarily essentially unique) functor $F\:\bbA\to\bbB$ such that $\Phi=\bbB(-,F-)$. This amounts to a 2-functor 
\begin{equation}\label{a2.0}
\Cat(\Q)\to\Map(\Dist(\Q))\:\Big(F\:\bbA\to\bbB\Big)\mapsto\Big(\bbB(-,F-)\:\bbA\dist\bbB\Big).
\end{equation}

We write $\Sym\Dist(\Q)$ for the full subquantaloid of $\Dist(\Q)$ determined by the symmetric $\Q$-categories. It is easily verified that the involution $f\mapsto f\o$ on the base quantaloid $\Q$ extends to the quantaloid $\Sym\Dist(\Q)$: explicitly, if $\Phi\:\bbA\dist\bbB$ is a distributor between symmetric $\Q$-categories, then so is $\Phi\o\:\bbB\dist\bbA$, with elements $\Phi\o(a,b):=\Phi(b,a)\o$. And if $F\:\bbA\to\bbB$ is a functor between symmetric $\Q$-categories, then the left adjoint distributor represented by $F$ has the particular feature that it is a symmetric left adjoint in $\Sym\Dist(\Q)$. That is to say, the functor in \eqref{a2.0} restricts to the symmetric situation as
\begin{equation}\label{a2.0.0}
\Sym\Cat(\Q)\to\Sym\Map(\Sym\Dist(\Q))\:\Big(F\:\bbA\to\bbB\Big)\mapsto\Big(\bbB(-,F-)\:\bbA\dist\bbB\Big),
\end{equation}
obviously giving a commutative diagram
$$\begin{array}{c}
\xymatrix@=8ex{
\Cat(\Q)\ar[r] & \Map(\Dist(\Q)) \\
\Sym\Cat(\Q)\ar[r]\ar[u]^{\mathrm{incl.}} & \Sym\Map(\Sym\Dist(\Q))\ar[u]_{\mathrm{incl.}}}
\end{array}$$

The full embedding $\Sym\Cat(\Q)\hookrightarrow\Cat(\Q)$ has a right adjoint functor:
$$\Sym\Cat(\Q)\xymatrix@=8ex{\ar@{}[r]|{\perp}\ar@<1mm>@/^2mm/[r]^{\mathrm{incl.}} & \ar@<1mm>@/^2mm/[l]^{(-)\s}}\Cat(\Q).$$
This {\em symmetrisation} sends a $\Q$-category $\bbA$ to the symmetric $\Q$-category $\bbA\s$ whose objects (and types) are those of $\bbA$, but for any two objects $x,y$ the hom-arrow is $\bbA\s(y,x):=\bbA(y,x)\wedge\bbA(x,y)\o$. A functor $F\:\bbA\to\bbB$ is sent to $F\s\:\bbA\s\to\bbB\s\:a\mapsto Fa$. It is a result of [Heymans and Stubbe, 2011] that the inclusion $\Sym\Map(\Sym\Dist(\Q))\to\Map(\Dist(\Q))$ admits a right adjoint that makes the diagram
$$\xymatrix@=8ex{
\Cat(\Q)\ar[r]\ar[d]_{(-)\s} & \Map(\Dist(\Q))\ar@{.>}[d] \\
\Sym\Cat(\Q)\ar[r] & \Sym\Map(\Sym\Dist(\Q))}$$
commute if and only if, for each family $(f_i\:X\to X_i, g_i\:X_i\to X)_{i\in I}$ of morphisms in $\Q$,
$$\left.\begin{array}{c}
\forall j,k\in I:\ f_k\circ g_j\circ f_j\leq f_k \\[1ex]
\forall j,k\in I:\ g_j\circ f_j\circ g_k\leq g_k \\[1ex]
1_X\leq\displaystyle\bigvee_{i\in I}g_i\circ f_i
\end{array}\right\}\Longrightarrow\
1_X\leq\bigvee_{i\in I}(g_i\wedge f_i\o)\circ(g_i\o\wedge f_i).$$
Such an involutive quantaloid $\Q$ is said to be {\em Cauchy-bilateral}. We will encounter examples of Cauchy-bilateral quantaloids further on in this paper.

\subsection*{Presheaves, Cauchy-completion and symmetric-completion}

A {\em (contravariant) presheaf} on $\bbA$ is a distributor into $\bbA$ whose domain is a one-object category with an identity hom-arrow. Writing $*_X$ for the one-object $\Q$-category whose single object $*$ has type $X\in\Q_0$ and whose single hom-arrow is the identity $1_X$, a presheaf is then typically written as $\phi\:*_X\dist\bbA$. The set of presheaves on $\bbA$ is written $\P(\bbA)$: it is a $\Q$-category when we define that $t(\phi\:*_X\dist\bbA)=X$ and $\P(\bbA)(\psi,\phi)=[\psi,\phi]$ (this being a lifting in the quantaloid $\Dist(\Q)$, i.e.\ the value at $\phi$ of the right adjoint to composition with $\psi$). The {\em Yoneda embedding} of $\bbA$ into $\P(\bbA)$ is the fully faithful functor of $\Q$-enriched categories $Y_{\bbA}\:\bbA\to\P(\bbA)$ that sends $a\in\bbA$ to the {\em representable presheaf} $\bbA(-,a)\:*_{ta}\dist\bbA$. In fact, this procedure extends to a functor $\P:\Cat(\Q)\to\Cat(\Q)$, which is the {\em free cocompletion KZ-doctrine} on the category of $\Q$-categories. (A {\em covariant presheaf} on $\bbA$ is a distributor $\phi\:\bbA\dist *_X$; they are not of much importance in this paper.)

A $\Q$-category $\bbA$ is said to be {\em Cauchy complete} when each left adjoint distributor with codomain $\bbA$ is represented by a functor [Lawvere, 1973], that is, when for each $\Q$-category $\bbB$ the functor in \eqref{a2.0} determines an equivalence of ordered sets
$$\Cat(\Q)(\bbB,\bbA)\simeq\Map(\Dist(\Q))(\bbB,\bbA).$$ 
This clearly implies that the functor in \eqref{a2.0} restricts to a biequivalence of locally ordered 2-categories between $\Cat\cc(\Q)$, the full subcategory of $\Cat(\Q)$ determined by the Cauchy complete $\Q$-categories, and $\Map(\Dist(\Q))$. Moreover, the full inclusion of $\Cat\cc(\Q)$ in $\Cat(\Q)$ admits a left adjoint:
\begin{equation}\label{z1}
\Cat\cc(\Q)\xymatrix@=8ex{\ar@{}[r]|{\perp}\ar@<-1mm>@/_2mm/[r]_{\mathrm{full\ incl.}} & \ar@<-1mm>@/_2mm/[l]_{(-)\cc}}\Cat(\Q).
\end{equation}
That is to say, each $\Q$-category $\bbA$ has a {\em Cauchy completion} $\bbA\cc$: it is the full subcategory of the presheaf category $\P(\bbA)$ whose objects are the left adjoint presheaves on $\bbA$. The Yoneda embedding $Y_{\bbA}\:\bbA\to\P(\bbA)$ factors through $\bbA\cc$, and the distributor induced by $Y_{\bbA}\:\bbA\to\bbA\cc$ turns out to be an isomorphism in $\Dist(\Q)$. Therefore the quantaloid $\Dist(\Q)$ is equivalent to its full subquantaloid $\Dist\cc(\Q)$ whose objects are the Cauchy complete $\Q$-categories. As a result, there is an equivalence of locally ordered 2-categories
\begin{equation}\label{xxx1}
\Cat\cc(\Q)\simeq\Map(\Dist\cc(\Q))\simeq\Map(\Dist(\Q)).
\end{equation}

The Cauchy completion can of course be applied to a symmetric $\Q$-category (assuming that $\Q$ is involutive), but the resulting Cauchy complete category need not be symmetric anymore: the functor $(-)\cc\:\Cat(\Q)\to\Cat(\Q)$ does not restrict to $\Sym\Cat(\Q)$ in general. However, its very definition suggests the following modification [Heymans and Stubbe, 2011]: a symmetric $\Q$-category $\bbA$ is {\em symmetrically complete} if, for any symmetric $\Q$-category $\bbB$, the functor in \eqref{a2.0.0} determines an equivalence of symmetrically ordered sets
$$\Sym\Cat(\Q)(\bbB,\bbA)\simeq\Sym\Map(\Sym\Dist(\Q))(\bbB,\bbA).$$
This implies that the functor in \eqref{a2.0.0} restricts to a biequivalence between $\Sym\Cat\sc(\Q)$, the full subcategory of $\Sym\Cat(\Q)$ of its symmetrically complete objects, and $\Sym\Map(\Sym\Dist(\Q))$. Moreover, the full inclusion of $\Sym\Cat\sc(\Q)$ in $\Sym\Cat(\Q)$ admits a left adjoint:
\begin{equation}\label{z3}
\Sym\Cat\sc(\Q)\xymatrix@=8ex{\ar@{}[r]|{\perp}\ar@<-1mm>@/_2mm/[r]_{\mathrm{full\ incl.}} & \ar@<-1mm>@/_2mm/[l]_{(-)\sc}}\Sym\Cat(\Q).
\end{equation}
Explicitly, for a symmetric $\Q$-category $\bbA$, its {\em symmetric completion} $\bbA\sc$ is the full subcategory of the Cauch completion $\bbA\cc$ (and thus also a full subcategory of the presheaf category $\P(\bbA)$) determined by the {\em symmetric} left adjoint presheaves. For similar reasons as above, there is an equivalence of involutive quantaloids between $\Sym\Dist(\Q)$ and its full subquantaloid $\Sym\Dist\sc(\Q)$ of symmetrically complete $\Q$-categories, and therefore also an equivalence of categories
\begin{equation}\label{xxx2}
\Sym\Cat\sc(\Q)\simeq\Sym\Map(\Sym\Dist\sc(\Q))\simeq\Sym\Map(\Sym\Dist(\Q)).
\end{equation}

Importantly, a result of [Heymans and Stubbe, 2011] says that, if $\Q$ is a Cauchy-bilateral quantaloid, then the symmetric-completion and the Cauchy-completion of any symmetric $\Q$-category coincide, and the symmetrisation of a Cauchy complete $\Q$-category is symmetrically complete. In fact, there is a distributive law of the monad $(-)\cc\:\Cat(\Q)\to\Cat(\Q)$ over the comonad $(-)\s\:\Cat(\Q)\to\Cat(\Q)$. This means in particular that there is a commutative diagram of adjunctions as follows:
$$\xymatrix@=8ex{
\Cat\cc(\Q)\ar@<-1mm>@/_2mm/[d]_{(-)\s}\ar@{=}[r] & \Map(\Dist\cc(\Q))\ar@<-1mm>@/_2mm/[d]_{(-)\s}\ar@{=}[r] & \Map(\Dist(\Q))\ar@<-1mm>@/_2mm/[d]_{(-)\s} \\
\Sym\Cat\sc(\Q)\ar@<-1mm>@/_2mm/[u]_{\mbox{incl.}}\ar@{}[u]|{\vdash}\ar@{=}[r] & \Sym\Map(\Sym\Dist\sc(\Q))\ar@<-1mm>@/_2mm/[u]_{\mbox{incl.}}\ar@{}[u]|{\vdash}\ar@{=}[r] & \Sym\Map(\Sym\Dist(\Q))\ar@<-1mm>@/_2mm/[u]_{\mbox{incl.}}\ar@{}[u]|{\vdash}}$$
The equal signs in this diagram are the equivalences of \eqref{xxx1} and \eqref{xxx2}; the bottom row is fully included in the top row, and can be obtained from it by `symmetrisation'.

\subsection*{Universal constructions}

An idempotent in a quantaloid $\Q$ is, of course, an endomorphism $e\:A\to A$ such that $e^2=e$. Such an idempotent is said to split in $\Q$ when there exists a diagram
\begin{equation}\label{abc}
\xymatrix@=8ex{
A\ar@(l,u)^{qp=e}\ar@<1ex>[r]^p & B\ar@(u,r)^{1_B=pq}\ar@<1ex>[l]^q}
\end{equation}
in $\Q$. If $\E$ is a class of idempotents in a quantaloid $\Q$, then we write $\Q_{\E}$ for the quantaloid obtained by splitting the idempotents in $\E$. An explicit description goes as follows: the objects of $\Q_{\E}$ are the elements of $\E$, and $\Q_{\E}(e,f)=\{x\:A\to B\mid f\circ x=x=x\circ e\}$ whenever $e\:A\to A$ and $f\:B\to B$ are in $\E$. Composition and local suprema in $\Q_{\E}$ are as in $\Q$, but the identity on an idempotent $e$ is, obviously, $e\:e\to e$ itself. If all identities in $\Q$ are in $\E$, then there is a fully faithful homomorphism of quantaloids
$$I\:\Q\to\Q_{\E}\:\Big(x\:A\to B\Big)\mapsto\Big(x\:1_A\to 1_B\Big)$$ 
which is the universal splitting in $\Q$ of idempotents in $\E$. Spelled out, this means that if $F\:\Q\to\R$ is a homomorphism of quantaloids, and the images of all idempotents in $\E$ split in $\R$, then there is an essentially unique homomorphism $\overline{F}\:\Q_{\E}\to\R$ such that $\overline{F}\circ I=F$. Moreover, if $F$ is fully faithful then so is $\overline{F}$.

When $\Q$ is an involutive quantaloid, then we say that an idempotent $e\:A\to A$ in $\Q$ is {\em symmetric} when $e\o=e$. It is straightforward that $\Q_{\E}$ is then involutive too: the involute of $x\in\Q_{\E}(e,f)$ is computed as in $\Q$, for the symmetry of $e$ and $f$ make sure that $x\o\in\Q_{\E}(f,e)$. As before, it is surely the case that, whenever all identities in $\Q$ are in $\E$, the involutive quantaloid $\Q_{\E}$ has a universal property for the splitting of idempotents. Noting however that $I\:\Q\to\Q_{\E}$ preserves the involution, we can point out a slightly more subtle universal property. Say that the splitting in the diagram in \eqref{abc} is {\em symmetric} when $q=p\o$. If $F\:\Q\to\R$ is a homomorphism of involutive quantaloids and the images of all idempotents in $\E$ split symmetrically in $\R$, then there is an essentially unique homomorphism $\overline{F}\:\Q_{\E}\to\R$ of involutive quantaloids such that $\overline{F}\circ I=F$; in other words, if $F$ preserves the involution then so does $\overline{F}$. And again, if $F$ is fully faithful then so is $\overline{F}$. 

In any quantaloid $\Q$, products and sums are the same thing, so they are usually referred to as {\em direct sums}. We write $A=\oplus_{i\in I}A_i$ for the direct sum of a family $(A_i)_i$ of objects of $\Q$, with injections $s_i\:A_i\to A$ and projections $p_i\:A\to A_i$; in fact, for $A$ to be the direct sum of the $(A_i)_i$, it is a necessary and sufficient condition that $p_i\circ s_j=\delta_{ij}$ and $\bigvee_{i}s_i\circ p_i=1_A$. In these equations, $\delta_{ij}\:A_j\to A_i$ is the ``Kronecker delta'': it is the identity morphism when $i=j$ and the zero morphism otherwise. The universal direct sum completion of a small quantaloid $\Q$ exists, and can explicitly be described as the quantaloid $\Matr(\Q)$ of {\em matrices} over $\Q$. An object in $\Matr(\Q)$ is a $\Q$-typed set, i.e.\ a set $A$ together with a type function $t\:A\to\Q_0$, and a morphism between two such $\Q$-typed sets is a {\em matrix} $M\:A\to B$, i.e.\ a family $M(b,a)\:ta\to tb$ of morphisms in $\Q$, one for each $(a,b)\in A\times B$. Of course, matrices can be composed: for $M\:A\to B$ and $N\:B\to C$ we have $N\circ M\:A\to C$ with elements 
$$(N\circ M)(c,a):=\bigvee_{b\in B}N(c,b)\circ M(b,a).$$
The identity on a $\Q$-typed set $A$ is the matrix $\Delta_A\:A\to A$ all of whose elements are ``Kronecker deltas''. With elementwise supremum, this makes $\Matr(\Q)$ a quantaloid; and whenever $\Q$ is involutive, so is $\Matr(\Q)$ (for elementwise involution). There is a fully faithful homomorphism
$$J\:\Q\to\Matr(\Q)\:\Big(f\:X\to Y\Big)\mapsto\Big((f):\{X\}\to\{Y\}\Big)$$
sending a morphism to the matrix between singletons in the obvious way (which preserves the involution on $\Q$ whenever there is one), which is the universal direct sum completion of $\Q$. 

Any $\Q$-typed set $A$ determines a $\Q$-category $\bbA$ by putting $\bbA_0=A$ and $\bbA(a',a)=\Delta_A(a',a)$: this is precisely a {\em discrete} $\Q$-category in the sense that the hom-arrow between two different objects is a zero morphism and every endo-hom-arrow is an identity morphism. A matrix between $\Q$-typed sets is easily seen to be precisely a distributor between discrete $\Q$-categories, so the quantaloid $\Matr(\Q)$ is precisely the full subquantaloid of $\Dist(\Q)$ of discrete $\Q$-categories. A discrete $\Q$-category is obviously symmetric, so whenever $\Q$ is an involutive quantaloid, $\Matr(\Q)$ can also be considered as full involutive subcategory of $\Sym\Dist(\Q)$. Furthermore, a monad in $\Matr(\Q)$ is exactly a $\Q$-category, and (assuming that $\Q$ is involutive) a symmetric monad is a symmetric $\Q$-category. In other words, both $\Dist(\Q)$ and $\Sym\Dist(\Q)$ can be constructed from $\Matr(\Q)$ by splitting a particular class of idempotents:
\begin{itemize}
\item $\Dist(\Q)=\Matr(\Q)_{\E}$ for $\E$ the class of monads in $\Matr(\Q)$,
\item $\Sym\Dist(\Q)=\Matr(\Q)_{\E}$ for $\E$ the class of symmetric monads in $\Matr(\Q)$.
\end{itemize}
Composing the various universal constructions we thus find how $\Dist(\Q)$ and $\Sym\Dist(\Q)$ can be considered as completions of $\Q$ itself. 

For any involutive quantaloid $\Q'$ it is a matter of fact that the process of splitting {\em all} monads in $\Q'$ can be broken down in two steps: first split all {\em symmetric} monads in $\Q'$, then split all {\em anti-symmetric} monads in the thusly obtained quantaloid. (A monad $m\:X\to X$ in an involutive quantaloid is said to be anti-symmetric when $m\wedge m\o=1_X$.) Applying this to $\Q'=\Matr(\Q)$ for a small involutive quantaloid $\Q$, this exhibits how $\Dist(\Q)$ is also a completion of $\Sym\Dist(\Q)$.

All this goes to show that both $\Dist(\Q)$ and $\Sym\Dist(\Q)$ lead a ``double life''. On the one hand, they are concretely constructed quantaloids: their objects are (symmetric) $\Q$-categories, and their morphisms are distributors. This makes it possible to compute with {\em individual objects and morphisms} of $\Dist(\Q)$ (or $\Sym\Dist(\Q)$). But on the other hand, $\Dist(\Q)$ and $\Sym\Dist(\Q)$ are universal constructions on $\Q$: first add all direct sums to $\Q$, then split either all monads or only the symmetric ones. These universal properties thus say something about {\em the collection of all objects and morphisms} of $\Dist(\Q)$ or $\Sym\Dist(\Q)$. The first approach is clearly rooted in the theory of quantaloid-enriched categories, whereas the second approach is close in spirit to allegory theory. Indeed, quoting P. Johnstone [2002, p.\ 138], ``many allegories of interest may be generated by idempotent-splitting processes from quite small full sub-allegories''. Of course, $\Dist(\Q)$ or $\Sym\Dist(\Q)$ need not be allegories (neither of them is necessarily modular, see further), but they are both generated by universal processes from a quite small full sub-quantaloid, namely from $\Q$ itself.

\subsection*{Orders and sheaves over a base quantaloid}

We now have everything ready to state the central definitions with which we shall work in this paper. First we recall a definition first given in [Stubbe, 2005b]:
\begin{definition}\label{C3.1}
Given a small quantaloid $\Q$ and a set $\E$ of idempotents in $\Q$, we define
$$\Ord(\Q,\E):=\Cat\cc(\Q_{\E})\mbox{ \ and \ }\Idl(\Q,\E):=\Dist\cc(\Q_{\E})$$
for, respectively, the locally ordered 2-category of {\em $(\Q,\E)$-orders and order functions}, and the quantaloid of {\em $(\Q,\E)$-orders and ideal relations}. If $\E$ is taken to be the set of {\em all} idempotents in $\Q$, then we write $\Q\si$ instead of $\Q_{\E}$, $\Ord(\Q)$ instead of $\Ord(\Q,\E)$, and $\Idl(\Q)$ instead of $\Idl(\Q,\E)$; we then simply speak of {\em $\Q$-orders} (and order functions and ideal relations).
\end{definition}

Next we present a new definition, intended as ``symmetric'' version of the previous definition. Because the term ``symmetric $\Q$-order'' is technically inadequate (it suggests a $\Q$-order with a symmetric hom, {\em quod non}), and the term ``$\Q$-set'' already means something related-but-different in the literature (see e.g.\ [Higgs, 1973; Fourman and Scott, 1979; Borceux, 1994; Mulvey and Nawaz, 1995; Gylys, 2001; Johnstone, 2002; and others]), we opt to speak of ``$\Q$-sheaves'':
\begin{definition}\label{C3}
Given a small involutive quantaloid $\Q$ and a set $\E$ of symmetric idempotents in $\Q$, we define
$$\Set(\Q,\E):=\Sym\Cat\sc(\Q_{\E})\mbox{ \ and \ }\Rel(\Q,\E):=\Sym\Dist\sc(\Q_{\E})$$
for, respectively, the category of {\em $(\Q,\E)$-sheaves and functions}, and the quantaloid of {\em $(\Q,\E)$-sheaves and relations}. If $\E$ is taken to be the set of {\em all} symmetric idempotents in $\Q$, then we write $\Q\ssi$ instead of $\Q_{\E}$, $\Set(\Q)$ instead of $\Set(\Q,\E)$, and $\Rel(\Q)$ instead of $\Rel(\Q,\E)$; we then simply speak of {\em $\Q$-sheaves} (and functions and relations).
\end{definition}
We shall explain at the end of Section \ref{D} how, for so-called small quantaloids of closed cribles, the symmetry condition in the above definition is in fact equivalent to an appropriate discreteness condition.

From the general theory on (symmetric) $\Q$-categories that we explained in the previous subsections, we can now conclude that:
\begin{proposition}\label{C4}
For any small quantaloid $\Q$ and any set $\E$ of idempotents in $\Q$, there is a biequivalence of locally ordered 2-categories
$$\Ord(\Q,\E)\xymatrix{\ar[r]^{\sim}&}\Map(\Idl(\Q,\E))\:\Big(F\:\bbA\to\bbB\Big)\mapsto\Big(\bbB(-,F-)\:\bbA\dist\bbB\Big).$$
For any small involutive quantaloid $\Q$ and any set $\E$ of symmetric idempotents in $\Q$, there is an equivalence of categories 
$$\Set(\Q,\E)\xymatrix{\ar[r]^{\sim}&}\Sym\Map(\Rel(\Q,\E))\:\Big(F\:\bbA\to\bbB\Big)\mapsto\Big(\bbB(-,F-)\:\bbA\dist\bbB\Big).$$
If $\Q$ is an involutive quantaloid and $\E$ a set of symmetric idempotents such that $\Q_{\E}$ is Cauchy-bilateral, then both squares in
$$\xymatrix@=8ex{
\Ord(\Q,\E)\ar@<-1mm>@/_2mm/[d]_{(-)\s}\ar[r]^(0.4){\sim} & \Map(\Idl(\Q,\E))\ar@<-1mm>@/_2mm/[d]_{(-)\s} \\
\Set(\Q,\E)\ar[r]^(0.4){\sim}\ar@<-1mm>@/_2mm/[u]_{\mbox{incl.}}\ar@{}[u]|{\vdash} & \Sym\Map(\Rel(\Q,\E))\ar@<-1mm>@/_2mm/[u]_{\mbox{incl.}}\ar@{}[u]|{\vdash}}$$
commute, and the bottom row is obtained by ``symmetrising'' the top row.
\end{proposition}

Here is yet another result of the general theory of $\Q$-categories:
\begin{proposition}\label{C5}
For any small (resp.\ involutive) quantaloid $\Q$ and any set $\E$ of (resp.\ symmetric) idempotents in $\Q$, there is an equivalence of (resp.\ involutive) quantaloids 
$$\Idl(\Q,\E)\simeq\Dist(\Q_{\E})\mbox{, resp.\ }\Rel(\Q,\E)\simeq\Sym\Dist(\Q_{\E}).$$
\end{proposition}
This proposition explains an important subtlety: each $(\Q,\E)$-order (or $(\Q,\E)$-sheaf) is {\em Morita equivalent} to a (symmetric) $\Q_{\E}$-category. This fact has often been used (implicitly) to forget about Cauchy completeness altogether: several definitions of ``sheaf on an involutive quantaloid'' that can be found in the literature, amount (in one form or another) to stating that a sheaf is a symmetric category, and a morphism of sheaves is a left adjoint distributor. (An example that springs to mind, is the formalism of {\em projection matrices}, on which we shall comment in more detail in Section \ref{E}.) However, we have deliberately opted to include the requirement of Cauchy (or symmetric) completeness in the definition of ``sheaf'' on a quantaloid $\Q$, for it expresses precisely the ``gluing condition'' that one expects of such a notion (as well illustrated by [Walters, 1981]). But of course it comes in handy that, modulo Morita equivalence, this completeness can be swiped under the carpet.

The whole of Section \ref{D} is devoted to showing that the topos of sheaves on a site $(\C,J)$ is equivalent to $\Set(\Q)$ when $\Q=\R(\C,J)$ is the small quantaloid of closed cribles in $(\C,J)$.

\section{Sheaves on a site}\label{D}

For any small involutive quantaloid $\Q$ we stated in Definition \ref{C3} and Proposition \ref{C4} that
$$\Set(\Q):=\Sym\Cat\sc(\Q\ssi)\simeq\Sym\Map(\Sym\Dist(\Q\ssi)).$$ 
Walters [1982] showed that, for a small site $(\C,J)$,    
$$\Sh(\C,J)\simeq\Sym\Cat\cc(\R(\C,J))\simeq\Map(\Sym\Dist(\R(\C,J))),$$
where $\R(\C,J)$ is the so-called {\em small quantaloid of closed cribles} (which Walters originally called the {\em bicategory of relations}) constructed from $(\C,J)$. In this section we shall show that sheaves on a small site $(\C,J)$ (in the topos-theoretic sense) correspond with sheaves on the small involutive quantaloid $\R(\C,J)$ (in the sense of our Definition \ref{C3}). 

More precisely, we shall prove that, if $\Q=\R(\C,J)$ is a small quantaloid of closed cribles, then for any set $\E$ of symmetric idempotents in $\Q$ containing the identities, $\Sym\Dist(\Q)$ and $\Sym\Dist(\Q_{\E})$ are equivalent modular quantaloids; and because each left adjoint in a modular quantaloid is necessarily a symmetric left adjoint, it follows that $\Set(\Q,\E)$ is equivalent to $\Map(\Sym\Dist(\Q))$, which in turn is equivalent to $\Sh(\C,J)$ by Walters' [1982] result. To give our proof, we shall use the axiomatic description of $\R(\C,J)$ due to [Heymans and Stubbe, 2012], for it allows us to prove our claim via elementary computations in involutive quantaloids, much in the line of Freyd and Scedrov's [1990] work on allegories (see also [Johnstone, 2002]). In the next subsection we recall the necessary results from our earlier work.

\subsection*{Axioms for a small quantaloid of closed cribles}

First we recall some definitions:
\begin{definition}\label{D2}
A quantaloid $\Q$ is:
\begin{enumerate}
\item\label{D2.1} {\em locally localic} if, for all objects $X$ and $Y$, $\Q(X,Y)$ is a locale,
\item\label{D2.2} {\em map-discrete} if, for any left adjoints $f\:X\to Y$ and $g\:X\to Y$ in $\Q$, $f\leq g$ implies $f=g$,
\item\label{D2.3} {\em weakly tabular} if, for every $q\:X\to Y$ in $\Q$, 
$$q=\bigvee\left\{fg^*\bigmid (f,g)\:X\span Y\mbox{ is a span of left adjoints such that }fg^*\leq q\right\},$$
\item\label{D2.4} {\em map-tabular} if for every $q\:X\to Y$ in $\Q$ there is a span $(f,g)\:X\span Y$ of left adjoints in $\Q$ such that $fg^*=q$ and $f^*f\wedge g^*g=1_{\dom(f)}$,
\item\label{D2.5} {\em weakly modular} if, for every pair of spans of left adjoints in $\Q$, say $(f,g)\:X\span Y$ and $(m,n)\:X\span Y$, we have $fg^*\wedge mn^*\leq f(g^*n\wedge f^*m)n^*$,
\item\label{D2.6} {\em tabular} if it is involutive and if for every $q\:X\to Y$ in $\Q$ there exists a span $(f,g)\:X\span Y$  of left adjoints in $\Q$ such that $fg\o=q$ and $f\o f\wedge g\o g=1_{\dom(f)}$,
\item\label{D2.7} {\em modular} if it is involutive and if for any $f\:X\to Y$, $g\:Y\to Z$ and $h\:X\to Z$ in $\Q$ we have $gf\wedge h\leq g(f\wedge g\o h)$ (or equivalently, $gf\wedge h\leq (g\wedge hf\o)f$).
\end{enumerate} 
\end{definition}
The notions of modularity\footnote{
In fact, J. Riguet [1948, p.\ 120] discovered much earlier what he called the {\em Dedekind formula} (in French, {\it la relation de Dedekind}) for relations between sets: if $R\subseteq E\times F$, $S\subseteq F\times G$ and $T\subseteq E\times G$, then $SR\cap T\subseteq(S\cap TR\o)(R\cap S\o T)$ (where $R\o$ is the opposite relation of $R$, etc.). Whereas it is obvious that the Dedekind formula implies the modular law, it is not difficult to see that the converse holds too: $SR\cap T=(SR\cap T)\cap(SR\cap T)\subseteq(SR\cap T)\cap(S\cap TR\o)R\subseteq(S\cap TR\o)((S\cap TR\o)\o(SR\cap T)\cap R)\subseteq(S\cap TR\o)(S\o T\cap R)$. All this can, of course, be done in any involutive locally ordered 2-category, and indeed Riguet certainly understood that the importance of the Dedekind formula went beyond the calculus of relations: he explains that the term {\em relation de Dedekind} was deliberately so chosen because ``{\it elle contient comme cas particulier la relation entre idéaux dans un anneau découverte par Dedekind}''.} and tabularity are cited from Freyd and Scedrov [1990] who give them in the context of allegories\footnote{Freyd and Scedrov [1990] define an {\em allegory} $\A$ to be a locally posetal 2-category, equipped with an involution $\A\to\A\op\:f\mapsto f\o$ (which fixes the objects, reverses the arrows, and preserves the local order), in which the modular law holds. Johnstone [2002] calls an allegory {\em geometric} when its hom-posets are complete lattices and composition distributes over arbitrary suprema. Thus, a geometric allegory is exactly the same thing as a modular quantaloid.}. Weak modularity, weak tabularity and map-tabularity were introduced in [Heymans and Stubbe, 2012] with the specific aim to axiomatise small quantaloids of closed cribles.

There are many useful relations between several of these notions; we recall some of these in the next lemma.
\begin{lemma}\label{D3}
\begin{enumerate}
\item In any modular quantaloid, all left adjoints are symmetric left adjoints.
\item Any modular quantaloid is map-discrete.
\item Any locally localic and modular quantaloid is Cauchy-bilateral.
\item A small quantaloid $\Q$ is weakly tabular if and only if $\Dist(\Q)$ is map-tabular.
\item A small quantaloid $\Q$ is locally localic and modular if and only if $\Matr(\Q)$ is modular.
\end{enumerate}
\end{lemma} 
The first two statements in the above lemma appear in [Freyd and Scedrov, 1990; Johnstone, 2002], the third is quoted from [Heymans and Stubbe, 2011], and the two other statements come from [Heymans and Stubbe, 2012]. Also the following result appears in the latter reference.
\begin{theorem}\label{D1}
For a small quantaloid $\Q$, the following conditions are equivalent:
\begin{enumerate}
\item $\Q$ is a {\em small quantaloid of closed cribles}, i.e.\ there exists a small site such that $\Q\simeq\R(\C,J)$, 
\item putting, for $X\in\Map(\Q)$,
$$J(X):=\Big\{S\mbox{ is a sieve on }X\bigmid 1_X=\bigvee_{s\in S}ss^*\Big\}$$
defines a Grothendieck topology $J$ on $\Map(\Q)$ for which $\Q\cong\R(\Map(\Q),J)$,
\item $\Q$ is locally localic, map-discrete, weakly tabular and weakly modular. 
\end{enumerate}
In this case, $\Q$ carries an involution, sending $q\:Y\to X$ to
$$q\o:=\bigvee\left\{gf^*\bigmid (f,g)\:Y\span X\mbox{ is a span of left adjoints such that }fg^*\leq q\right\},$$
which makes $\Q$ also modular.
\end{theorem}

\subsection*{Splitting symmetric idempotents}

In this subsection we study the properties of the involutive quantaloid $\Sym\Dist(\Q)$ when $\Q$ is a small quantaloid of closed cribles.

First we point out two useful conditions to determine whether a (small or large) involutive quantaloid $\Q'$ has symmetric splittings for its symmetric idempotents. The first lemma can be found in [Freyd and Scedrov, 1990, 2.162; Johnstone, 2002, Lemma A3.3.3] and the second is a corollary of [Freyd and Scedrov, 1990, 2.166 and 2.169; Johnstone, 2002, A3.3.6 and A3.3.12]. For completeness' sake we shall give proofs here too, specifically adapted to the situation at hand.
\begin{lemma}\label{D8.1.0}
If $\Q'$ is a modular quantaloid, then each splitting of a symmetric idempotent is necessarily a symmetric splitting.
\end{lemma}
\proof
First observe that for any $f\:X\to Y$ in a modular $\Q'$ we always have $f\leq ff\o f$:  because $f=1_Yf\wedge f\leq(1_Y\wedge ff\o)f\leq ff\o f$. Now suppose that $e\:A\to A$, $p\:A\to B$ and $q\:B\to A$ satisfy $e=e^2=e\o=qp$ and $pq=1_B$ in $\Q'$. Then it follows that $q\o=pqq\o\leq pp\o pqq\o =pp\o q\o=p(qp)\o=pe\o=pe=p$, and similarly $p\o\leq q$.
\endofproof
\begin{lemma}\label{D8}
If $\Q'$ is a modular and tabular quantaloid in which all symmetric monads\footnote{In the context of allegories, Freyd and Scedrov [1990] use the term {\em equivalence relation} (and Johnstone [2002] speaks simply of an {\em equivalence}) for what we call a symmetric monad; when it splits, then it does so symmetrically (because an allegory is modular), and they say that the equivalence relation is {\em effective}. If all equivalence relations in an allegory split, they say that the allegory is effective.} split, then all symmetric idempotents in $\Q'$ have a (necessarily symmetric) splitting.
\end{lemma}
\proof
Let $e\:A\to A$ be a symmetric idempotent in $\Q'$: we shall exhibit a splitting. To that end, first consider a tabulation $(f,g)$ of $e\wedge 1_A$:
$$\xymatrix@R=8ex@C=4ex{
 & B\ar[dl]_g\ar[dr]^f \\
A\ar[rr]_{1_A\wedge e} & & A}$$
Thus, $f$ and $g$ are left adjoints in $\Q'$ such that $fg\o=1_A\wedge e$ and $g\o g\wedge f\o f=1_B$. Because $\Q'$ is modular, we know moreover that $f\dashv f\o$ and $g\dashv g\o$. It is useful to point out that $g\leq ef$ and $f\leq eg$ follow from these assumptions, and that, in turn, this implies that $g\leq eg$ and $eg=ef$.

Now define $t:=(eg)\o(eg)=g\o eg\:B\to B$. Then clearly $t\o=t$ holds; it is furthermore easy to check that $tt=g\o egg\o eg\leq g\o e1_Aeg=g\o eg=t$; and $t=(eg)\o eg\geq g\o g\geq 1_B$ follows from inequalities pointed out above. In sum, this says that $t\:B\to B$ is a symmetric monad. By assumption we can split $t$: there is a diagram
$$\xymatrix@=8ex{
B^t\ar@(l,u)^{1_{B^t}}\ar@<1ex>[r]^{h\o} & B\ar@<1ex>[l]^{h}\ar@(u,r)^t}$$
such that $t=h\o h$ and $hh\o=1_{B^t}$ (where, again by modularity, $h\dashv h\o$). 

Next, consider the diagram
$$\xymatrix@=8ex{
B\ar@(l,u)^t\ar@<1ex>[r]^{eg} & A\ar@<1ex>[l]^{(eg)\o}\ar@(u,r)^e}$$
in which, by definition of $t$, we have $t=(eg)\o(eg)$. Using modularity of $\Q'$ and the tabulation $(f,g)$ of $1_A\wedge e$, we can compute that $e=e1_Ae\wedge e\leq e(1_A\wedge e\o ee\o)e=e(1_A\wedge e)e=e(fg\o)e=(ef)(eg)\o=(eg)(eg)\o$. But $(eg)(eg)\o=egg\o e\leq ee=e$ follows immediately from $g\dashv g\o$, hence we obtain $e=(eg)(eg)\o$. 

Composing these two diagrams produces a splitting in $\Q'$ for the symmetric idempotent $e\:A\to A$, as required.
\endofproof
For any small involutive quantaloid $\Q$, $\Sym\Dist(\Q)$ is a quantaloid in which all symmetric monads split: simply because it is the universal splitting of symmetric monads in $\Matr(\Q)$. Below we shall furthermore prove that, whenever $\Q$ is a small quantaloid of closed cribles, $\Sym\Dist(\Q)$ is modular and tabular too.
\begin{lemma}\label{D4}
If $\Q$ is a locally localic quantaloid and $\E$ is a collection of idempotents in $\Q$, then $\Q_{\E}$ is locally localic too.
\end{lemma}
\proof
If $p\:L\to L$ is an idempotent sup-morphism on a complete lattice, then $p(L)\subseteq L$ is a complete lattice too, with the same suprema as in $L$, but with $p(x)\wedge'p(y):=p(p(x)\wedge p(y))$ as binary infimum and $p(\top)$ as emtpy infimum (i.e.\ top element). A simple computation shows that, if $L$ is a locale, then so is $p(L)$. This applies to $e_2\circ-\circ e_1\:\Q(X,Y)\to\Q(X,Y)$ for idempotents $e_1\:X\to X$ and $e_2\:Y\to Y$ in $\E$, to show that $\Q_{\E}(e_1,e_2)$ is a locale whenever $\Q(X,Y)$ is; hence $\Q_{\E}$ is locally localic whenever $\Q$ is.
\endofproof
\begin{lemma}\label{D5}
If $\Q$ is a modular quantaloid and $\E$ is a collection of symmetric idempotents in $\Q$, then $\Q_{\E}$ is modular too.
\end{lemma}
\proof
Local suprema, composition and involution in $\Q_{\E}$ are the same as in $\Q$. As pointed out in the above proof, the infimum of $f,g\:e_1\to e_2$ in $\Q_{\E}$ is $f\wedge'g:=e_2(f\wedge g)e_1$, but thanks to the modular law it is easily seen that 
$$f\wedge g=e_2f\wedge ge_1\leq e_2(fe_1\o\wedge e_2\o g)e_1=e_2(f\wedge g)e_1,$$
whereas $e_2(f\wedge g)e_1\leq f\wedge g$ is always valid, hence in this case the local binary infima in $\Q_{\E}$ are the same as in $\Q$. Thus it follows that $\Q_{\E}$ is modular whenever $\Q$ is.
\endofproof
\begin{proposition}\label{D6}
If $\Q$ is a small, locally localic, modular quantaloid, then $\Sym\Dist(\Q)$ is modular.
\end{proposition}
\proof
$\Matr(\Q)$ is modular by Lemma \ref{D3}, so $\Sym\Dist(\Q)=(\Matr(\Q))_{\E}$, with $\E$ the collection of symmetric monads in $\Matr(\Q)$, is modular too by Lemma \ref{D5}.
\endofproof
\begin{proposition}\label{D7}
If $\Q$ is a small, weakly tabular, Cauchy-bilateral quantaloid, then $\Sym\Dist(\Q)$ is tabular.
\end{proposition}
\proof
From [Heymans and Stubbe, 2012, Proposition 3.5] we recall that a small quantaloid $\Q$ is weakly tabular if and only if $\Dist(\Q)$ is map-tabular. The proof for the necessity goes as follows: Suppose that $\Phi\:\bbA\dist\bbB$ is a distributor. We can assume without loss of generality that $\bbA$ and $\bbB$ are Cauchy complete, because every $\Q$-category is isomorphic to its Cauchy completion in $\Dist(\Q)$. Now define the $\Q$-category $\bbR$ to be the full subcategory of $\bbA\times\bbB$ whose objects are those $(a,b)\in\bbA\times\bbB$ for which $1_{ta}\leq\Phi(a,b)$, and write $T$ (resp.\ $S$) for the composition of the inclusion $\bbR\hookrightarrow\bbA\times\bbB$ with the projection of $\bbA\times\bbB$ onto $\bbA$ (resp.\ onto $\bbB$). By construction we then have $\bbB(S-,S-)\wedge\bbA(T-,T-)=\bbR$; and, relying on the weak tabularity of $\Q$ and the Cauchy completeness of $\bbA$ and $\bbB$, a lenghty computation shows that $\Phi=\bbA(-,T-)\tensor\bbB(S-,-)$. That is to say, the left adjoints $\bbA(-,T-)\:\bbR\dist\bbA$ and $\bbB(-,S-)\:\bbR\dist\bbB$ in $\Dist(\Q)$ provide for a map-tabulation of $\Phi\:\bbA\dist\bbB$. 

We now modify this proof to suit our needs. For any $\Phi\:\bbB\dist\bbA$ in $\Sym\Dist(\Q)$ we must find $\Sigma\:\bbR\dist\bbA$ and $\Theta\:\bbR\dist\bbB$ in $\Sym\Dist(\Q)$ such that $\Sigma\tensor\Theta\o=\Phi$ and $\Sigma\o\tensor\Sigma\ \wedge\ \Theta\o\tensor\Theta=\bbR$. If $\Q$ is Cauchy-biateral then the Cauchy completion of a symmetric $\Q$-category is again symmetric, hence any symmetric $\Q$-category is isomorphic to its Cauchy completion in $\Sym\Dist(\Q)$ (and not merely in $\Dist(\Q)$). Therefore we may still suppose that $\bbA$ and $\bbB$ are Cauchy complete. Referring to the above, the category $\bbR$ is clearly symmetric whenever $\bbA$ and $\bbB$ are, and the left adjoint distributors represented by the functors $S\:\bbR\to\bbB$ and $T\:\bbR\to\bbA$ are evidently symmetric left adjoints. Thus the result follows.
\endofproof
In view of Theorem \ref{D1} we may now conclude from the above:
\begin{theorem}\label{D7.0}
If $\Q$ is a small quantaloid of closed cribles, then $\Sym\Dist(\Q)$ is a modular and tabular quantaloid in which all symmetric idempotents split symmetrically.
\end{theorem}

\subsection*{Change of base}

This subsection is devoted to the proof of the fact that, when $\Q$ is a small quantaloid of closed cribles and $\E$ is a class of symmetric idempotents in $\Q$ containing all identities, then also $\Q_{\E}$ is a small quantaloid of closed cribles, and the involutive quantaloids $\Sym\Dist(\Q)$ and $\Sym\Dist(\Q_{\E})$ are equivalent. To tackle this problem, we study the ``change of base'' homomorphism from $\Sym\Dist(\Q)$ to $\Sym\Dist(\Q_{\E})$ which is determined by the universal property of splitting symmetric idempotents. Let us first recall the appropriate terminology.

Let $\Q$ and $\Q'$ be small involutive quantaloids and $F\:\Q\to\Q'$ be a homomorphism that preserves the involution. It is easily seen that a symmetric $\Q$-category $\bbA$ determines a symmetric $\Q'$-category $F\bbA$ by putting:
\begin{itemize}
\item objects: $(F\bbA)_0=\bbA_0$ with types $t_{F\bbA}a=F(ta)$ in $\Q'_0$,
\item hom-arrows: $(F\bbA)(a',a)=F(\bbA(a',a))$ for all objects $a,a'$.
\end{itemize}
Similarly for distributors, and $F$ so determines a homomorphism $\overline{F}\:\Sym\Dist(\Q)\to\Sym\Dist(\Q')$ of involutive quantaloids that makes the diagram
\begin{equation}\label{D7.2.0}
\xymatrix@=8ex{
\Q\ar[r]^F\ar@{^{(}->}[d]_I & \Q'\ar@{^{(}->}[d]^{I'} \\
\Sym\Dist(\Q)\ar[r]_{\overline{F}} & \Sym\Dist(\Q')}
\end{equation}
commute: $\overline{F}$ is the {\em change of base} homomorphism induced by $F$. (Of course, $I$ denotes the canonical inclusion of $\Q$ in $\Sym\Dist(\Q)$, and similarly for $I'$.)

Now we recall a necessary and sufficient condition for the ``change of base'' induced by some $F\:\Q\to\Q'$ to be an equivalence. As it is straightforward to verify that $F\:\Q\to\Q'$ is fully faithful if and only if the change of base $\overline{F}$ is fully faithful, we need to take a closer look at the essential surjectivity on objects of $\overline{F}$. 
\begin{lemma}\label{D7.2}
Let $F\:\Q\to\Q'$ be a homomorphism of small involutive quantaloids. The change of base $\overline{F}\:\Sym\Dist(\Q)\to\Sym\Dist(\Q')$ is an equivalence of involutive quantaloids if and only if there exists a fully faithful homomorphism of involutive quantaloids $G\:\Q'\to\Sym\Dist(\Q)$ making the diagram below essentially commutative:
$$\xymatrix@=8ex{
\Q\ar[r]^F\ar@{^{(}->}[d]_I & \Q'\ar[dl]^G \\
\Sym\Dist(\Q)}$$
\end{lemma}
\proof
First suppose that $\overline{F}$ is an equivalence. Considering the commutative square in \eqref{D7.2.0}, the required fully faithful $G$ is obtained by composing $I'$ with the pseudo-inverse of $\overline{F}$.

Conversely, suppose that a fully faithful $G$ exists such that $G\circ F\cong I$. Because $G$ and $I$ are fully faithful, so is $F$, and thus also $\overline{F}$. Size issues apart, also $I$ and $G$ induce a change of base, and we end up with an essentially commutative diagram
$$\xymatrix@=8ex{
\Sym\Dist(\Q)\ar[r]^{\overline{F}}\ar[d]_{\overline{I}} & \Sym\Dist(\Q')\ar[dl]^{\overline{G}} \\
\Sym\Dist(\Sym\Dist(\Q))}$$
The homomorphisms $\overline{I}$, $\overline{F}$ and $\overline{G}$ are fully faithful, because $I$, $F$ and $G$ are. If we show that $\overline{I}$ is essentially surjective on objects, then it is an equivalence, and hence so is $\overline{F}$.

To see that $\overline{I}$ is indeed essentially surjective on objects, one can do as follows. Given $\bbC$ in $\Sym\Dist(\Sym\Dist(\Q))$, let us explicitly write the hom-arrow from an object $x\in\bbC$ to an object $y\in\bbC$ as $\Gamma_{y,x}\:\bbA_x\dist\bbA_y$; these morphisms in $\Sym\Dist(\Q)$ satisfy the conditions that make $\bbC$ a symmetric category: $\bbA_x\leq\Gamma_{x,x}$, $\bigvee_{y\in\bbC}\Gamma_{z,y}\tensor\Gamma_{y,x}\leq\Gamma_{z,x}$ and $\Gamma_{x,y}=\Gamma_{y,x}\o$ (for all $x,y\in\bbC$). With these data, we define a symmetric $\Q$-category $\bbA$ as follows:
\begin{itemize}
\item objects: $\bbA_0:=\biguplus_{x\in\bbC}\bbA_x$, with inherited types,
\item hom-arrows: for $u\in\bbA_x$ and $v\in\bbA_y$, $\bbA(v,u):=\Gamma_{y,x}(v,u)$.
\end{itemize}
Regarding $\bbA$ now as an object in $\Sym\Dist(\Sym\Dist(\Q))$, via the change of base $\overline{I}$, we further define a distributor $\Gamma\:\overline{I}(\bbA)\dist\bbC$ by:
\begin{itemize}
\item distributor-elements: for $u\in\bbA_x$ and $y\in\bbC$, $\Gamma(y,u):=\Gamma_{y,x}(-,u)$.
\end{itemize}
It is then a fact that $\Gamma\tensor\Gamma\o=\bbC$ and $\Gamma\o\tensor\Gamma=\overline{I}(\bbA)$. All verifications are long but straightforward computations.
\endofproof
In the exact same situation as in the above lemma, we can sometimes say more:
\begin{lemma}\label{D8.2.1}
In the situation of Lemma \ref{D7.2}, if $G$ is fully faithful and $\Sym\Dist(\Q)$ is modular and tabular then $\Q'$ is modular and weakly tabular.
\end{lemma}
\proof
Modularity of $\Q'$ follows straightforwardly from the modularity of $\Sym\Dist(\Q)$ and the fully faithful homomorphism $\Q'\to\Sym\Dist(\Q)$ of involutive quantaloids.

To deduce the weak tabularity of $\Q'$ from the tabularity and modularity of $\Sym\Dist(\Q)$, we first make a helpful observation. Given any $\Phi\:\bbA\dist\bbB$ in $\Sym\Dist(\Q)$, let
$$\xymatrix@R=8ex@C=4ex{
 & \bbC\ar[dl]|{\distsign}_{\Sigma}\ar[dr]|{\distsign}^{\Theta} \\
\bbA\ar[rr]|{\distsign}_{\Phi} & & \bbB}$$
be a tabulation; then, in particular, $\Phi=\Theta\tensor\Sigma\o$ and $\Sigma\dashv\Sigma\o$. Now consider the family
$$\Big(\bbC(-,c)\:*_{tc}\dist\bbC\Big)_{c\in\bbC}$$
of all representable presheaves on $\bbC$, each of which is a left adjoint in $\Sym\Dist(\Q)$. Precomposing both $\Sigma\:\bbC\dist\bbA$ and $\Theta\:\bbC\dist\bbB$ with these thus gives a family, indexed by the $c\in\bbC$, 
$$\xymatrix@R=8ex@C=4ex{
 & {*_{tc}}\ar[dl]|{\distsign}_{\Sigma\tensor\bbC(-,c)=\Sigma(-,c)}\ar[dr]|{\distsign}^{\Theta(-,c)=\Theta\tensor\bbC(-,c)} \\
\bbA & & \bbB}$$
of spans of left adjoints in $\Sym\Dist(\Q)$, whose domains are in the image of the canonical embedding $\Q\hookrightarrow\Sym\Dist(\Q)$, such that
$$\Phi=\bigvee_{c\in\bbC}\Theta(-,c)\tensor\Big(\Sigma(-,c)\Big)^*.$$
In particular, if $\Phi\:\bbA\dist\bbB$ is in the image of $G\:\Q'\to\Sym\Dist(\Q)$, then -- because the image of $\Q\hookrightarrow\Sym\Dist(\Q)$ is contained in the image of $G$ -- it admits a weak tabulation by spans of left adjoints in the image of $G$. By fully faithfulness of $G$, $\Q'$ is weakly tabular.
\endofproof
The above results apply in particular when $\Q$ is a small quantaloid of closed cribles and when we put $\Q'=\Q\ssi$: they show that splitting the symmetric idempotents in a small quantaloid of closed cribles is ``harmless'' for the theory of sheaves. In fact, instead of splitting {\em all} symmetric idempotents, we can choose to split only those in a class $\E$ of symmetric idempotents containing all identities.
\begin{theorem}\label{D8.4} If $\Q$ is a small quantaloid of closed cribles and $\E$ is a class of symmetric idempotents in $\Q$ containing all identities, then also $\Q_{\E}$ is a small quantaloid of closed cribles and the inclusion $\Q\hookrightarrow\Q_{\E}$ induces an equivalence $\Sym\Dist(\Q)\simeq\Sym\Dist(\Q_{\E})$ of involutive quantaloids.
\end{theorem}
\proof
If $\Q$ is a small quantaloid of closed cribles, then it is locally localic, hence so is $\Q_{\E}$, by Lemma \ref{D4}. The other results follow from the commutative diagram
$$\xymatrix@=8ex{
\Q\ar@{^{(}->}[r]\ar@{^{(}->}[d] & \Q_{\E}\ar[dl] \\
\Sym\Dist(\Q)
}$$
of fully faithful functors, and the fact that $\Sym\Dist(\Q)$ is modular and tabular.
\endofproof
Of course, taking $\E$ to be the class of all symmetric idempotents in $\Q$, we find that $\Q\ssi$ is a small quantaloid of closed cribles such that $\Sym\Dist(\Q)\simeq\Sym\Dist(\Q\ssi)$. But taking $\E$ to be the class of all symmetric monads in $\Q$, or the class of all symmetric comonads\footnote{In a modular quantaloid $\Q$, an arrow $c\:A\to A$ is a symmetric comonad if and only if it satisfies $c\leq 1_A$. Indeed, if $c\leq 1_A$ then $c\o\leq 1_A$ too, hence $c\leq cc\o c$ (cf.\ the proof of Lemma \ref{D8.1.0}) implies that $c\leq cc$ but also  that $c\leq c\o$, and by involution also $c\o\leq c$. Therefore, particularly in the context of allegories, the term {\em coreflexive} is often used.}, produces other important examples.

\subsection*{Walters' theorem revisited}

We now have everything ready to make the following extension to the result of [Walters, 1982]. As is customary, we write $\Rel(\topE)$ for the quantaloid of internal relations in a topos $\topE$.
The next theorem excludes all confusion with our earlier notation $\Rel(\Q)$.
\begin{theorem}\label{D9}
For any small site $(\C,J)$, any small quantaloid $\Q\simeq\R(\C,J)$ and any set $\E$ of symmetric idempotents in $\Q$ containing all identities, we have the following equivalences:
\begin{enumerate}
\item $\Set(\Q,\E)\simeq\Sym\Cat\cc(\Q)\simeq\Sh(\C,J)$,
\item $\Rel(\Q,\E)\simeq\Sym\Dist(\Q)\simeq\Rel(\Sh(\C,J))$.
\end{enumerate}
\end{theorem}
\proof
This proof relies on Walters' [1982, p.~101] theorem that the topos $\Sh(\C,J)$ is biequivalent to the bicategory $\Sym\Cat\cc(\R(\C,J))$ (Walters' insistence on the term {\em biequivalence} stresses the fact that a single morphism in the category $\Sh(\C,J)$ gets identified with an equivalence class of morphisms in the bicategory $\Sym\Cat\cc(\R(\C,J))$ whose homs are symmetric preorders), on Freyd and Scecrov's [1990, 2.148] theorem that any tabular allegory $\A$ is equivalent (as allegory) to the allegory $\Rel(\Map(\A))$ of internal relations in the regular category $\Map(\A)$, and on the particular properties of $\Sym\Dist(\Q)$, for $\Q$ a small quantaloid of closed cribles, that we summarised in Theorems \ref{D7.0} and \ref{D8.4}.

Because $\Q$ is a small quantaloid of closed relations, so is $\Q_{\E}$ (see Theorem \ref{D8.4}); one particular consequence is that, for symmetric categories enriched in either quantaloid, the symmetric completion coincides with the Cauchy completion (cf.\ Lemma \ref{D3} and Proposition \ref{C5}). Furthermore, again by Theorem \ref{D8.4}, $\Sym\Dist(\Q)$ is equivalent to $\Sym\Dist(\Q_{\E})$. All this justifies the following equivalences of involutive quantaloids:
$$\Rel(\Q,\E):=\Sym\Dist\sc(\Q_{\E})=\Sym\Dist\cc(\Q_{\E})\simeq\Sym\Dist(\Q_{\E})\simeq\Sym\Dist(\Q).$$
By Theorem \ref{D7.0} we know that $\Sym\Dist(\Q)$ is a modular quantaloid, hence so are its equivalents; all left adjoints in the above involutive quantaloids are therefore symmetric left adjoints, by Lemma \ref{D3}. Taking (symmetric) left adjoints therefore produces the following equivalences of categories (or rather, biequivalences of 2-categories which are locally symmetrically ordered):
$$\Set(\Q,\E):=\Cat\sc(\Q_{\E})\simeq\Sym\Map(\Sym\Dist\sc(\Q_{\E}))\simeq\Map(\Sym\Dist(\Q))\simeq\Sym\Cat\cc(\Q).$$
Invoking at this point Walters' theorem, this proves (1). But because the involutive quantaloid $\Sym\Dist(\Q)$ is not only modular but also tabular (see again Theorem \ref{D7.0}), Freyd and Scedrov's theorem proves it to be equivalent to the involutive quantaloid of internal relations in $\Map(\Sym\Dist(\Q))$, which in turn proves (2). 
\endofproof
The theorem above thus says two things about a small quantaloid of closed cribles $\Q$ and a set $\E$ of symmetric idempotents in $\Q$ containing all identities: firstly, that the category $\Set(\Q,\E):=\Sym\Cat\sc(\Q_{\E})$ is the category of sheaves on a site; secondly, that this category $\Set(\Q,\E)$ admits, up to equivalence, the simpler description $\Sym\Cat\cc(\Q)$. (And similar for $\Rel(\Q,\E)$.) Choosing $\E$ to be the set of all symmetric idempotents in $\Q$, we find:
\begin{corollary}
If $\Q$ is a small quantaloid of closed cribles, then $\Set(\Q)$ is a Grothendieck topos and $\Rel(\Q)$ is its category of relations. 
\end{corollary}

\subsection*{Symmetric vs.\ discrete}

In this subsection we wish to make a remark on the symmetry axiom that we used in Definition \ref{C3} of $\Q$-sheaves. In any locally ordered category $\K$, an object $D$ is said to be {\em discrete} when, for any other object $X\in\K$, the order $\K(X,D)$ is symmetric. In [Heymans and Stubbe, 2012] we showed that, for a Cauchy-bilateral quantaloid $\Q$, every symmetric and Cauchy complete $\Q$-category is a discrete object of $\Cat\cc(\Q)$. In general the converse need not hold, but:
\begin{proposition} If $\Q$ is a small quantaloid of closed cribles, then a Cauchy complete $\Q$-category is discrete in $\Cat\cc(\Q)$ if and only if it is symmetric.
\end{proposition}  
\proof
Suppose that $\bbA$ is a discrete object in $\Cat\cc(\Q)$;  we seek to prove that $\bbA(y,x)=\bbA(x,y)\o$ for any $x,y\in\bbA$. Relying in particular on the weak tabularity of $\Q$, it is sufficient to show that, for any span $(f,g)\:ty\span tx$ of left adjoints in $\Q$, 
$$fg^*\leq\bbA(x,y)\iff fg^*\leq\bbA(y,x)\o.$$
But, because $\bbA$ is Cauchy complete, for any such span $(f,g)$ we can consider the tensors $x\tensor f$ and $y\tensor g$ in $\bbA$, and writing $U=\dom(f)=\dom(g)$ we indeed have 
\begin{eqnarray*}
f\circ g^*\leq\bbA(x,y)
 & \iff & 1_U\leq\bbA(x\tensor f,y\tensor g) \\
 & \iff & 1_U\leq\bbA(y\tensor g,x\tensor f) \\
 & \iff & g\circ f^*\leq\bbA(y,x)\\
 & \iff & f\circ g^*\leq\bbA(y,x)\o
\end{eqnarray*}
where the second equivalence is due to the discreteness of $\bbA$ and the last equivalence holds because $f^*=f\o$ due to the modularity of $\Q$.
\endofproof
If $\Q$ is a small quantaloid of closed cribles, then so is $\Q\ssi$, and the Cauchy competion and symmetric completion of a symmetric $\Q\ssi$-enriched category coincide. Thus we find:
\begin{corollary} If $\Q$ is a small quantaloid of closed cribles, then $\Set(\Q)$ is the full subcategory of discrete objects of $\Ord(\Q)$ and $\Rel(\Q)$ is the full subquantaloid of discrete objects of $\Idl(\Q)$.
\end{corollary}
That is to say, whereas we defined the objects of $\Set(\Q)$ as the {\em symmetric} objects in $\Ord(\Q)$, we now find that they are exactly the {\em discrete} objects.
 
\section{Grothendieck quantaloids and quantales}\label{E}

In the previous section we showed that, for $\Q$ a small quantaloid of closed cribles, $\Set(\Q)$ is a Grothendieck topos and $\Rel(\Q)$ is its category of relations. Given this result, it is a natural to ask whether this is the case for other involutive quantaloids too; and if so, for which ones. Precisely, we wish to find necessary and sufficient conditions on $\Q$ for $\Rel(\Q)$ to be the category of relations in a topos.
\begin{definition}\label{10}
A small involutive quantaloid $\Q$ is called a {\em Grothendieck quantaloid} (if $\Q$ has only one object we speak of a {\em Grothendieck quantale}) if there exists a topos $\topE$ such that there is an equivalence $\Rel(\topE)\simeq\Rel(\Q)$ of involutive quantaloids.
\end{definition} 
A sufficient condition on $\Q$ is being a small quantaloid of closed cribles. On the other hand, the internal relations in a topos form a modular quantaloid, and $\Q$ is a full subquantaloid of $\Rel(\Q)$, so a necessary condition will be the modularity of $\Q$. To establish a precise necessary-and-sufficient condition, we first point out a connection with {\em projection matrices}.
\begin{definition}
If $\Q$ is a small involutive quantaloid, then $\Proj(\Q):=\Matr(\Q)\ssi$ is the involutive quantaloid of {\em projection matrices}\footnote{In [Freyd and Scedrov, 1990, 2.226], in the context of allegories rather than quantaloids, this construction is referred to as the {\em systemic completion}.}.
\end{definition}
Straightforwardly extending the terminology for quantales [Resende, 2012], we shall say that a quantaloid $\Q$ is {\em stably Gelfand} if it is an involutive quantaloid in which $ff\o f\leq f$ implies $f\leq ff\o f$ for any morphism $f\:X\to Y$. Any modular quantaloid is trivially stably Gelfand, as seen in the proof of Lemma \ref{D8.1.0}.
\begin{lemma}\label{11}
\begin{enumerate}
\item\label{ee1} For a stably Gelfand quantaloid $\Q$ there is an equivalence $\Proj(\Q)\simeq\Sym\Dist(\Q\ssi)\simeq\Rel(\Q)$ of involutive quantaloids. 
\item\label{ee2} If $\Q$ is a small quantaloid of closed cribles, then $\Proj(\Q)$ and $\Sym\Dist(\Q)$ are equivalent involutive quantaloids.
\item\label{ee3} A small involutive quantaloid $\Q$ is a Grothendieck quantaloid if and only if there exists a topos $\topE$ such that there is an equivalence $\Rel(\topE)\simeq\Proj(\Q)$ of involutive quantaloids.
\end{enumerate}
\end{lemma}
\proof
(\ref{ee1}) Let $P\:X\to X$ be a symmetric idempotent in $\Matr(\Q)$; that is to say, $X$ is a $\Q$-typed set, and $P$ is a collection of $\Q$-morphisms $P(x',x):tx\to tx'$, one for each $(x,x')\in X\times X$, such that 
$$\bigvee_{x''\in X}P(x',x'')\circ P(x'',x)=P(x',x)=P(x,x')\o\mbox{ for every $(x,x')\in X$.}$$
From this it is clear that $P(x,x')\circ P(x,x')\o\circ P(x,x')\leq P(x,x')$, so that by hypothesis  the converse inequality holds too. The computation
\begin{eqnarray*}
P(x,x')
 & \leq & P(x,x')\circ P(x,x')\o\circ P(x,x') \\
 & = & P(x,x')\circ P(x',x)\circ P(x,x') \\
 & \leq & P(x,x)\circ P(x,x') \\
 & \leq & P(x,x').
\end{eqnarray*}
then shows that $P(x,x')=P(x,x)\circ P(x,x')$; and similarly for $P(x,x')=P(x,x')\circ P(x',x')$. In other words, each $P(x,x)$ is an object of $\Q\ssi$, and each $P(x,x')$ is a morphism in $\Q\ssi$ from $P(x',x')$ to $P(x,x)$. As a consequence, we can define a symmetric $\Q\ssi$-category $\bbP$ whose $\Q\ssi$-typed object set is $X$ with types $tx:=P(x,x)$, and whose hom-arrows are $\bbP(x,x'):=P(x,x')$. Note that the $\Q\ssi$-category $\bbP$ is {\em normal} in the sense of [Stubbe, 2005b]: all of its endo-hom-arrows are identities. Furthermore, if both $P\:X\to X$ and $Q\:Y\to Y$ are projection matrices, and $M\:P\to Q$ is a morphism in $\Proj(\Q)$, i.e.\ a matrix $M\:X\to Y$ such that $Q\circ M=M=M\circ P$, then we can define a distributor $\Phi\:\bbP\dist\bbQ$ with elements $\Phi(y,x)=M(y,x)$. In fact, each distributor between $\bbP$ and $\bbQ$ arises in this way. In short, the correspondence $P\mapsto\bbP$ extends to an equivalence of involutive quantaloids between $\Proj(\Q)$ and the full involutive subquantaloid of $\Sym\Dist(\Q\ssi)$ of the {\em normal} symmetric $\Q\ssi$-categories (compare with [Stubbe, 2005b, Lemma 6.1]). But furthermore, a long but straightforward computation shows that each symmetric $\Q\ssi$-category is Morita equivalent with a normal symmetric $\Q\ssi$-category: so $\Sym\Dist(\Q\ssi)$ is equivalent to its full involutive subquantaloid of {\em normal} objects (compare with [Stubbe, 2005b, Lemma 6.2]). Taken together, all this proves that the correspondence $P\mapsto\bbP$ extends to an equivalence of involutive quantaloids between $\Proj(\Q)$ and $\Sym\Dist(\Q\ssi)$. Finally, by Proposition \ref{C5} the latter is furthermore equivalent to $\Rel(\Q):=\Sym\Dist\sc(\Q\ssi)$ (as involutive quantaloid).
 
(\ref{ee2}) Holds by Theorem \ref{D9}, taking $\E$ to be the set of all symmetric idempotents in $\Q$.

(\ref{ee3}) If $\Q$ is a Grothendieck quantaloid, then $\Rel(\Q)\simeq\Rel(\topE)$ for some topos $\topE$, so $\Q$ is modular because it is a full involutive subquantaloid of $\Rel(\Q)$. If, on the other hand, we assume that $\Proj(\Q)\simeq\Rel(\topE)$ for some topos $\topE$, then again $\Q$ is modular, now because it is a full involutive subquantaloid of $\Proj(\Q)$. In either case, $\Q$ is certainly stably Gelfand, so $\Proj(\Q)\simeq\Rel(\Q)$ by the first statement in this Lemma, which proves $\Proj(\Q)\simeq\Rel(\topE)\simeq\Rel(\Q)$ in either case.
\endofproof

We shall now recall some notions that [Freyd and Scedrov, 1990, 2.216(1), 2.225] introduced in the context of allegories, but that we adopt here for quantaloids. (In fact, the property that we call `weak semi-simplicity' was not given a name by Freyd and Scedrov [1990].)
\begin{definition}\label{12}
A morphism $q\:X\to Y$ in an involutive quantaloid $\Q$ is:
\begin{enumerate}
\item {\em simple} if $qq\o\leq 1_Y$,
\item {\em semi-simple} if there are simple morphisms $f$ and $g$ such that $q=fg\o$,
\item {\em weakly semi-simple} if $q=\bigvee\{fg\o\mid fg\o\leq q\mbox{ with $f$ and $g$ simple}\}$.
\end{enumerate}
And an involutive quantaloid $\Q$ is {\em (weakly) (semi-)simple} if each of its morphisms is.
\end{definition}
The next lemma can be found in [Freyd and Scedrov, 1990, 2.16(10)], but we spell out its proof for later reference.
\begin{lemma}\label{13}
A modular quantaloid $\Q$ is semi-simple if and only if $\Q\ssi$ is tabular.
\end{lemma}
\proof 
Suppose that $\Q\ssi$ is tabular. Given a morphism $q\:X\to Y$ in $\Q$, it can be included in $\Q\ssi$ as $q\:1_X\to 1_Y$, so let $q=fg\o$ be a tabulation in $\Q\ssi$. If $\Q$ is modular then so is $\Q\ssi$ (by Lemma \ref{D5}) so every left adjoint is a symmetric left adjoint. From this it is straightforward that both $f$ and $g$ are simple morphisms in $\Q$ such that $q=fg\o$.

Conversely, suppose first that $\Q$ is only semi-simple, and let $q\in\Q\ssi(r,p)$; that is to say, $r\:X\to X$ and $p\:Y\to Y$ are symmetric idempotents in $\Q$, and $q\:X\to Y$ is a morphism in $\Q$ satisfying $pq=q=qr$. Now let $b\:Z\to X$ and $a\:Z\to Y$ be simple morphisms in $\Q$ such that $ab\o=q$. It is then straightforward to check that $pa:1_Z\to p$ and $rb\:1_Z\to r$ are simple morphisms in $\Q\ssi$ such that $(pa)(rb)\o=q$. In other words, $\Q\ssi$ is semi-simple whenever $\Q$ is. Now the other hypothesis says that $\Q$ is also modular; by Lemma \ref{D5} we know that $\Q\ssi$ is modular too. It thus remains to prove that $\Q\ssi$ is tabular when it is semi-simple and modular. So again, let $q\in\Q\ssi(r,p)$ and suppose now that $x\:e\to p$ and $y\:e\to r$ are simple morphisms in $\Q\ssi$ such that $q=xy\o$. Simplicity of $x$ and $y$ makes $z:=x\o x\wedge y\o y\in\Q\ssi(e,e)$ a symmetric morphism satisfying $zz\leq z$; it follows from the modular law that it is therefore a symmetric idempotent, and furthermore that $xzy\o=xy\o$. Choosing a (necessarily symmetric) splitting of $z$ in $\Q\ssi$, say a $w\in\Q\ssi(f,e)$ such that $z=ww\o$ and $w\o w=f$, it is then tedious but routine to check that $x':=xw$ and $y':=yw$ are left adjoints in $\Q\ssi$ that tabulate $q$.
\endofproof
The proof of the previous lemma can be tweaked to obtain the following `weak' variant:
\begin{lemma}\label{13bis}
A modular quantaloid $\Q$ is weakly semi-simple if and only if $\Q\ssi$ is weakly tabular.
\end{lemma}
\proof 
This is a straightforward adaptation of the previous proof: instead of working with single pairs of simple morphisms we must work with families of pairs of simple morphisms.  

Suppose that $\Q\ssi$ is weakly tabular. A morphism $q\:X\to Y$ in $\Q$ can be viewed as a morphism $q\:1_X\to 1_Y$ in $\Q\ssi$, so consider its weak tabulation in $\Q\ssi$: 
$$q=\bigvee\{fg^*\mid \mbox{$f$ and $g$ are left adjoints in $\Q\ssi$ such that $fg^*\leq q$}\}.$$
If $\Q$ is modular then it follows (as in the previous proof) that all the $f$'s and $g$'s in the above expression are simple in $\Q$ and exhibit $q$'s weak semi-simplicity.

Conversely, suppose first that $\Q$ is weakly semi-simple. If $q\:r\to p$ is a morphism in $\Q\ssi$ (between symmetric idempotents $r\:X\to X$ and $p\:Y\to Y$, say) then at least we know that $q\:X\to Y$ is weakly semi-simple in $\Q$:
$$q=\{ab\o\mid\mbox{$a$ and $b$ are simple morphisms in $\Q$ such that $ab\o\leq q$}\}.$$
As in the previous proof, each such pair $(a,b)$ of simple morphisms in $\Q$ determines a pair $(pa,rb)$ of simple morphisms in $\Q\ssi$, and the lot of them exhibit $q$'s weak semi-simplicity in $\Q\ssi$. Thus $\Q\ssi$ is weakly semi-simple whenever $\Q$ is. Adding the hypothesis that $\Q$ is modular, we must prove that $\Q\ssi$ is in fact weakly tabular. So again, let $q\:r\to p$ be a morphism in $\Q\ssi$, and suppose now that
$$q=\bigvee\{xy\o\mid\mbox{$x$ and $y$ are simple morphisms in $\Q\ssi$ such that $xy\o\leq q$}\}.$$
Each of the pairs $(x,y)$ of simple morphisms in $\Q\ssi$ can be transformed, as in the previous proof, into a pair $(x',y')$ of left adjoint morphisms in $\Q\ssi$, and the lot of them provide for a weak tabulation of $q$. 
\endofproof

Much like Theorem \ref{D1} contains an axiomatic description of small quantaloids of closed cribles, we can now give an axiomatisation of Grothendieck quantaloids. In a sense, this is a refined analysis of the notion of `weak semi-simplicity'.
\begin{theorem}\label{14}
For a small involutive quantaloid $\Q$, the following are equivalent:
\begin{enumerate}
\item $\Q$ is weakly semi-simple,
\item $\Matr(\Q)$ is semi-simple,
\end{enumerate}
If $\Q$ is modular then this is also equivalent to:
\begin{enumerate}
\item[3.] $\Q\ssi$ is weakly tabular.
\end{enumerate}
If $\Q$ is modular and locally localic then this is also equivalent to:
\begin{enumerate}
\item[4.] $\Q\ssi$ is a small quantaloid of closed cribles,
\item[5.] $\Proj(\Q)$ is tabular,
\item[6.] there exists a small site $(\C,J)$ such that $\Rel(\Q)\simeq\Rel(\Sh(\C,J))$,
\item[7.] $\Q$ is a Grothendieck quantaloid.
\end{enumerate}
In fact, the small site $(\C,J)$ of which statement (6) speaks, is the site associated (as in Theorem \ref{D1}) with the small quantaloid of closed cribles $\Q\ssi$ of which statement (4) speaks.
\end{theorem}
\proof
($1\Rightarrow 2$) Let $M\:A\to B$ be a morphism in $\Matr(\Q)$: we must find semi-simple matrices $F\:C\to B$ and $G\:C\to A$ such that $M=FG\o$. Each element of $M$, that is, each $\Q$-arrow $M(b,a):ta\to tb$, is weakly semi-simple by assumption; thus
\begin{equation}\label{uvw}
M(b,a)=\bigvee\{fg\o\mid fg\o\leq M(b,a)\mbox{ with $f$ and $g$ simple}\}.
\end{equation}
For each $(a,b)\in A\times B$ we define the set 
$$C_{(a,b)}=\{(f,g)\mid fg\o\leq M(b,a)\mbox{ with $f$ and $g$ simple morphisms in $\Q$}\}$$
and furthermore we define $C$ to be the coproduct of the $C_{(a,b)}$'s. The constant functions $C_{(a,b)}\to A\:(f,g)\mapsto a$ and $C_{(a,b)}\to B\:(f,g)\mapsto b$ therefore uniquely define functions $\alpha\:C\to A$ and $\beta\:C\to A$; and putting the type of $(f,g)\in C$ to be the domain of $f$ (= the domain of $g$) makes $C$ an object of $\Matr(\Q)$. 

With the aid of the identity matrices $\Delta_A\:A\to A$ and $\Delta_B\:B\to B$ we now define two $\Q$-matrices, $F\:C\to B$ and $G\:C\to A$, to have as elements
$$F(b,(f,g))=\Delta_B(b,\beta(f,g))\circ f\mbox{ \ \ and \ \ }G(a,(f,g))=\Delta_A(a,\alpha(f,g))\circ g.$$
Simplicity of all $f$'s and $g$'s makes sure that $F$ and $G$ are simple matrices, and the formula $M=FG\o$ precisely coincides with \eqref{uvw}.

($2\Rightarrow 1$) Any $q\:X\to Y$ in $\Q$ may be viewed as a one-element matrix $(q)\:\{X\}\to\{Y\}$ between singletons (with obvious types). By hypothesis there are simple matrices $F\:C\to\{Y\}$ and $G\:C\to\{X\}$ such that $(q)=FG\o$. The simplicity of $F$ and $G$ implies that, for each $c\in C$, the morphisms $f_c:=F(Y,c)\:tc\to Y$ and $g_c:=G(c,X)\:X\to tc$ are simple morphisms in $\Q$; and $(q)=FG\o$ expresses precisely that $q=\bigvee_{c\in C}f_cg_c\o$, showing $q$ to be weakly semi-simple in $\Q$.

$(1\Leftrightarrow 3)$ This is the contents of Lemma \ref{13bis}

$(3\Leftrightarrow 4)$ If $\Q$ is a locally localic and modular quantaloid, then so is $\Q\ssi$ (by Lemmas \ref{D4} and \ref{D5}); in particular, $\Q\ssi$ is map-discrete and weakly modular too (see Lemma \ref{D3}). Thus $\Q\ssi$ is weakly tabular if and only if it is a small quantaloid of closed cribles (cf.\ Theorem \ref{D1}).

$(2\Leftrightarrow 5)$ From Lemmas \ref{D3} and \ref{D5} we know that $\Matr(\Q)$ is modular. Lemma \ref{13} does the rest, since $\Proj(\Q)=(\Matr(\Q))\ssi$.

$(4\Rightarrow 6)$ If $\Q\ssi\simeq\R(\C,J)$ for some small site $(\C,J)$, then $\Rel(\Q)\simeq\Rel(\Q\ssi)$ is equivalent to $\Rel(\Sh(\C,J))$ by Theorem \ref{D9}.

$(6\Rightarrow 7)$ Is evident.

$(7\Rightarrow 5)$ If $\Q$ is a Grothendieck quantaloid, then $\Proj(\Q)$ -- which by Lemma \ref{11} is equivalent to $\Rel(\Q)$ -- is equivalent to the allegory of internal relations in a topos, so it is most certainly tabular (see e.g.\ [Freyd and Scedrov, 1990, 2.142]).
\endofproof

\subsection*{Quantaloids vs.\ quantales}

The theorem above thus says that a Grothendieck quantaloid is precisely a modular, locally localic and weakly semi-simple quantaloid. There is an easier criterion than weak semi-simplicity when dealing with Grothendieck quantales rather than quantaloids. Using the term {\em quantal frame} to mean a quantale whose underlying sup-lattice is a locale [Resende, 2007], we can state it as:
\begin{theorem}\label{16}
A Grothendieck quantale is a modular quantal frame with a weakly semi-simple top (i.e.\ $\top=\bigvee\{f g\o\mid f,g\mbox{ simple}\}$).
\end{theorem}
\proof One implication is trivial. For the other, let $q\in Q$; we must show that it is weakly semi-simple. For any two simple elements of $Q$, $f$ and $g$ say, the modular law and the simplicity of $f$ and $g$ allow us to compute that
$$q\wedge fg\o\leq f(f\o qg\wedge 1)g\o\leq(ff\o qgg\o)\wedge f1g\o\leq q\wedge fg\o.$$
The element $h:=f(f\o qg\wedge 1)$ is simple, because it is smaller than the simple element $f$. In other words, this shows that, for any pair $(f,g)$ of simple elements, there exists a simple element $h$ such that $q\wedge fg\o=hg\o$. Using the remaining hypotheses, we can thus compute that
\begin{eqnarray*}
q & = & q\wedge\top \\
 & = & q\wedge\bigvee\{fg\o\mid\mbox{$f,g$ simple}\} \\
 & = & \bigvee\{q\wedge fg\o\mid\mbox{$f,g$ simple}\} \\
 & = & \bigvee\{hg\o\mid\mbox{$hg\o\leq q$ with $h$ and $g$ simple}\},
\end{eqnarray*}
so $Q$ is indeed weakly semi-simple.
\endofproof
As an application of the ``change of base'' principles that we developed in Section \ref{D}, we shall now show how every Grothendieck topos is equivalent to a category of $Q$-sheaves, with $Q$ a Grothendieck quantale.

First recall that two small quantaloids $\Q$ and $\R$ are said to be {\em Morita-equivalent} when the (large) quantaloids of modules $[\Q\op,\Sup]$ and $[\R\op,\Sup]$ are equivalent. B. Mesablishvili [2004] proved that for any small quantaloid $\Q$ there is a Morita-equivalent quantale $\Q\m$; he uses abstract $\V$-category theoretic arguments to prove his claim. Unraveling his arguments, we can give an explicit construction of $\Q\m$: it is $\Matr(\Q)(\Q_0,\Q_0)$, the quantale of endo-matrices with elements in $\Q$ on the $\Q$-typed set of objects of $\Q$ (where, of course, the type of an object $X\in\Q$ is $X$).

Given a morphism $f\:A\to B$ in a small quantaloid $\Q$, we shall write $M_f\in\Q\m$ for the matrix all of whose elements are zero, except for the element indexed by $(A,B)\in\Q_0\times\Q_0$, which is equal to $f$. The function $f\mapsto M_f$ is easily seen to preserve composition and suprema (but evidently not the identities, so it is not a quantaloid homomorphism). However, if $\E$ is a class of idempotents in $\Q\m$ containing all of $\{M_{1_A}\mid A\in \Q_0\}$, and we split these idempotents  in $\Q\m$, then we obtain a homomorphism
$$M\:\Q\to(\Q\m)_{\E}\:(f\:A\to B)\mapsto(M_f\:M_{1_A}\to M_{1_B})$$
which is easily seen to be fully faithful and injective on objects. If $\Q$ is a small involutive quantaloid, then it is straightforward to define an involution on the quantale $\Q\m$ as well, which the function $f\mapsto M_f$ preserves. If all elements of $\E$ are symmetric (which is automatic for the $M_{1_A}$), then the above homomorphism is not only fully faithful and injective on objects, but also preserves the involution.

Furthermore, $\Q\m$ is, by definition, a full subquantaloid of $\Matr(\Q)$, which in turn is a full subquantaloid of $\Sym\Dist(\Q)$; let us write the full inclusion as $J\:\Q\m\to\Sym\Dist(\Q)$. In case symmetric idempotents split symmetrically in $\Sym\Dist(\Q)$, there is a fully faithful homomorphism $J'\:(\Q\m)_{\E}\to\Sym\Dist(\Q)$ of involutive quantaloids. This is in particular the case when $\Q$ is a small quantaloid of closed cribles, which leads us to:
\begin{proposition}\label{D9.0}
If $\Q$ is a small quantaloid of closed cribles, $\Q\m$ is its Morita-equivalent quantale and $\E$ is a class of symmetric idempotents in $\Q\m$ containing all of $\{M_{1_A}\mid A\in \Q_0\}$, then also $(\Q\m)_{\E}$ is a small quantaloid of closed cribles and the inclusion $\Q\hookrightarrow(\Q\m)_{\E}$ induces an equivalence $\Sym\Dist(\Q)\to\Sym\Dist((\Q\m)_{\E})$ of involutive quantaloids.
\end{proposition}
\proof
If $\Q$ is a small quantaloid of closed cribles, then it is in particular locally localic and modular. Hence $\Matr(\Q)$ is locally localic, implying that $\Q\m$ is locally localic (as a one-object quantaloid), and therefore also $(\Q\m)_{\E}$ is locally localic. Moreover, it is straightforward to compute that the diagram
$$\xymatrix@=8ex{
\Q\ar@{^{(}->}[r]^M\ar@{^{(}->}[d]_I & (\Q\m)_{\E}\ar[dl]^{J'} \\
\Sym\Dist(\Q)
}$$
of involutive quantaloids and homomorphisms that preserve the involution commutes up to natural isomorphism. Because $J'$ is fully faithful, the results in Lemmas \ref{D7.2} and \ref{D8.2.1} apply, and prove the proposition.
\endofproof
If $\Q$ is a small quantaloid of closed cribles, then $\Sym\Dist(\Q\ssi)\simeq\Sym\Dist(\Q)$ by Proposition \ref{D8.4}, which is further equivalent to $\Sym\Dist((\Q\m)\ssi)$ by Proposition \ref{D9.0}. This produces the following:
\begin{corollary}\label{D9.1}
If $\Q$ is a small quantaloid of closed cribles, then
\begin{enumerate}
\item $\Set(\Q)\simeq\Set(\Q\m)$,
\item $\Rel(\Q)\simeq\Rel(\Q\m)$.
\end{enumerate}
This implies that $\Q\m$ is a Grothendieck quantale.
\end{corollary}
This result says in particular that {\em any Grothendieck topos can equivalently be described as a category of sets with an equality relation taking truth-values in a Grothendieck quantale.}

\subsection*{Examples}

We end this paper with some examples, the first two of which clearly illustrate the difference between `small quantaloids of closed cribles' and `Grothendieck quantaloids'.
\begin{example}[Closed cribles]\label{16.-1}
As remarked before, each small quantaloid $\Q$ of closed cribles is a Grothendieck quantaloid, and $\Set(\Q)$ is equivalent to the topos of sheaves on the site canonically associated with $\Q$. 
\end{example}
\begin{example}[Locales]\label{16.0}
A locale $(L,\bigvee,\wedge,\top)$ with its trivial involution is a Grothendieck quantale, but it is not a small quantaloid of closed cribles (because it is not weakly tabular). Upon splitting the (symmetric) idempotents in $L$ one obtains a small quantaloid of closed cribles; the site associated with the latter (as in Theorem \ref{D1}) is exactly the {\em canonical site} $(L,J)$ (for which $(x_i)_i\in J(x)$ if and only if $\bigvee_ix_i=x$). Thus $\Set(L)$ - in the sense of Definition \ref{C3} -- is equivalent to the ``usual'' topos of sheaves on $L$.
\end{example}

Our next example is somewhat more involved. First we must recall from [Stubbe, 2005b] that the 2-category $\Ord(\Q)$ of Definition \ref{C3.1} is equivalent to the 2-category ${\sf TRSCat}\cc(\Q)$ of ``Cauchy complete totally regular $\Q$-semicategories and totally regular semifunctors''; and in [Heymans and Stubbe, 2009a] it is shown to be further equivalent to the 2-category $\Map(\Mod_{\sf lpg}(\Q))$ of ``locally principally generated $\Q$-modules'' and left adjoint module morphisms. It is not difficult to deduce, from the symmetrisation of $\Q$-orders {\it qua} $\Q\ssi$-enriched categories that we proposed in Definition \ref{C3}, the appropriate symmetrisations of $\Q$-semicategories and of $\Q$-modules, thus producing as many different but equivalent descriptions of $\Q$-sheaves. In fact, in [Heymans and Stubbe, 2009b] we already studied the symmetric variant of locally principally generated $\Q$-modules, albeit only for involutive quantales (and not quantaloids): the so-called ``locally principally symmetric'' objects in $\Mod_{\sf lpg}(Q)$ form the subcategory $\Mod_{\sf lpg,lps}(Q)$. In [Heymans and Stubbe, 2009b, Example 3.7(4)] we showed that, for {\it any} involutive quantale $Q$, the involutive quantaloid $\Proj(Q)$ is equivalent to the involutive quantaloid $\Hilb(Q)$ of so-called ``$Q$-modules with Hilbert structure'' (and module morphisms between them). (The proof also appears in [Resende, 2012, Lemma 4.26, Theorem 4.29].) And we furthermore proved that, when $Q$ is a {\em modular quantal frame}, then $\Hilb(Q)$ is further equivalent to the involutive quantaloid $\Mod_{\sf lpg,lps}(Q)$ [Heymans and Stubbe, 2009b, Theorems 3.6 and 4.1]. Theorem \ref{14} says in particular that a Grothendieck quantale is necessarily a modular quantal frame, so together with Lemma \ref{11} this shows that in this case all of the involutive quantaloids $\Rel(Q)$, $\Proj(Q)$, $\Hilb(Q)$ and $\Mod_{\sf lpg,lps}(Q)$ are equivalent. Taking left adjoints in either of these therefore produces equivalent Grothendieck toposes
$$\Set(Q)\simeq\Map(\Rel(Q))\simeq\Map(\Proj(Q))\simeq\Map(\Hilb(Q))\simeq\Map(\Mod_{\sf lpg,lps}(Q)).$$
\begin{example}[Inverse quantal frames]\label{16.1}
An {\em inverse quantal frame} $Q$ is a modular quantal frame such that $\top=\bigvee\{p\in Q\mid p\o p\vee pp\o\leq 1\}$. It follows trivially from Theorem \ref{16} that inverse quantal frames are Grothendieck quantales. There is a correspondence up to isomorphism between inverse quantal frames and \'etale groupoids [Resende, 2007]: for every étale groupoid $G$ there is an inverse quantal frame $Q=\O(G)$ (as locale it is the object of morphisms of G, and its quantale multiplication stems from the composition law in $G$); and for every inverse quantal frame $Q$ there is an étale groupoid $G$ such that $Q\cong\O(G)$. Moreover, [Resende, 2012, p.~62--65] proves that $\Map(\Hilb(\O(G)))$ is equivalent to the classifying topos $BG$ of the étale groupoid $G$. Consequently,
$$\Set(\O(G))\!\simeq\!\Map(\Rel(\O(G)))\!\simeq\!\Map(\Proj(\O(G)))\!\simeq\!\Map(\Hilb(\O(G)))\!\simeq\!\Map(\Mod_{\sf lpg,lps}(\O(G)))$$ 
are all equivalent descriptions of the topos $BG$ in terms of ``sheaves on an involutive quantale''. 
\end{example}

\end{document}